\documentclass[12pt,a4paper]{article}

\usepackage[reqno]{amsmath}
\usepackage{amssymb,amsthm,upref,amscd}
\usepackage[T1]{fontenc}
\usepackage{times}
\usepackage{color}
\usepackage{array,multirow,makecell}

\newtheorem{theorem}{Theorem}[section]
\newtheorem{corollary}[theorem]{Corollary}

\newtheorem{lemma}[theorem]{Lemma}
\newtheorem{proposition}[theorem]{Proposition}
\newtheorem{remark}[theorem]{Remark}

\theoremstyle{definition} \theoremstyle{remark}
\numberwithin{equation}{section}

\begin{document}

\title{\textbf{Threshold for Existence, Non-existence and Multiplicity of positive solutions with prescribed mass for an NLS with a pure power nonlinearity in the exterior of a ball \ \\
}}
\author{Linjie Song$^{\mathrm{a,b,}}$\thanks{%
		Linjie Song is supported by CEMS. Email: songlinjie18@mails.ucas.edu.cn.}\  ,	
			Hichem Hajaiej$^{\mathrm{c}}$\thanks{%
				Email: hhajaie@calstatela.edu}\ \\
			{\small $^{\mathrm{a}}$Institute of Mathematics, AMSS, Academia Sinica,
				Beijing 100190, China}\\
			{\small $^{\mathrm{b}}$University of Chinese Academy of Science,
				Beijing 100049, China}\\
						{\small $^{\mathrm{c}}$ Department of Mathematics, California State University at Los Angeles,}\\
			{\small Los Angeles, CA 90032, USA}
		}
\date{}
\maketitle

\begin{abstract}

We obtain threshold results for the existence, non-existence and multiplicity of normalized solutions for semi-linear elliptic equations in the exterior of a ball. To the best of our knowledge, it is the first result in the literature addressing this problem.  In particular, we show that the prescribed mass can affect the number of normalized solutions and has a stabilizing effect in the mass supercritical case. Furthermore, in the threshold we find a new exponent $p = 6$ when $N = 2$, which does not seem to have played a role for this equation in the past. Moreover, our findings are "quite surprising" and completely different from the results obtained on the entire space and on balls. We will also show that the nature of the domain is crucial for the existence and stability of standing waves. As a foretaste, it is well-known that in the supercritical case these waves are unstable in $\mathbb{R}^N.$ In this paper, we will show that in the exterior domain they are strongly stable.

\bigskip



\noindent\textbf{Data availability statement:} The manuscript has no associate data.

\end{abstract}

\medskip

\section{Introduction}

In this paper, we are interested in the study of the standing wave solutions of the nonlinear Schr\"{o}dinger equation (NLS)
\begin{equation} \label{eq1.3}
\left\{
\begin{array}{cc}
i\partial_{t} \Psi + \Delta \Psi + |\Psi|^{p-2}\Psi = 0, (t,x) \in \mathbb{R} \times D, \\
\Psi(0,x) = \Psi_0(x),
\end{array}
\right.
\end{equation}
where $D = \{x \in \mathbb{R}^{N}: |x| > 1\}$, is the exterior of the unit ball in $\mathbb{R}^{N}$, $N \geq 2$, $2 < p < 2^{\ast},$ $2^{\ast} := 2N/(N-2)$ if $N \geq 3$ and $2^{\ast} := +\infty$ if $N = 2$.

A standing wave of \eqref{eq1.3} is a solution of the form $\Psi(t,x) = e^{i\lambda t}u(x)$ where the real-valued function $u$ solves $\Delta u - \lambda u + |u|^{p-2}u = 0$ in $D$, $u = 0$ on $\partial D$. We are interested in investigating the existence/non-existence and multiplicity of positive solutions for the following semi-linear elliptic equation
\begin{equation} \label{eq1.1}
\left\{
\begin{array}{cc}
-\Delta u + \lambda u = u^{p-1}, x \in D, \\
u > 0, u = 0 \ on \ \partial D, u(x) \rightarrow 0 \ as \ |x| \rightarrow +\infty,
\end{array}
\right.
\end{equation}
with $D = \{x \in \mathbb{R}^{N}: |x| > 1\}$, $N \geq 2$, $2 < p < 2^{\ast}$.

For \eqref{eq1.3}, it is well-known that we have the conservation of the energy $\mathcal{E}(\Psi)$ and the mass $\mathcal{M}(\Psi)$ given by:
$$
\mathcal{E}(\Psi) := \frac{1}{2}\int_{D}|\nabla \Psi|^{2}dx - \frac{1}{p}\int_{D}|\Psi|^{p}dx,
$$
$$
\mathcal{M}(\Psi) := \int_{D}|\Psi|^{2}dx.
$$
Therefore, a crucial step in studying the orbital stability of standing waves is to first establish the existence of the positive solutions of \eqref{eq1.1} satisfying:
\begin{equation} \label{eq1.2}
\int_{D}u^{2} = c,
\end{equation}
where $c > 0$ is prescribed. An exponent that has played an important role in the analysis of normalized solutions in $\mathbb{R}^N$ and balls is $p_{c} = 2+4/N.$ It is called in the literature the mass critical exponent. It separates two completely different regimes: the case $2 < p < 2 + 4/N$ is called mass subcritical and the case of $2 + 4/N < p < 2^{\ast}$ is called mass supercritical. There are many results dealing with the normalized solutions on $\mathbb{R}^{N}$, see e.g. \cite{HH, Stu1, Stu2, Stu3} for the mass subcritical case and \cite{ BN, BV, Jean} for the mass supercritical case. As for the mass supercritical case with non-autonomous nonlinearities, we have shown some results in two recent works \cite{Song, Song3}. However, all the techniques developed in these works do not apply to an exterior domain. To address this important case, we will develop a new method. Our main result reads as follows:

\begin{theorem}
\label{thm1.2} $(i)$ If $N \geq 3,$ and $2 + 4/N \leq p < 2^{\ast}$  or if $N = 2,$ and $4 < p < 6$, there exists $\eta_{1} > 0$ such that:\\
- For any $c > \eta_{1}$, (\ref{eq1.1}) admits at least two solutions $u_{\lambda}$ and $u_{\widetilde{\lambda}}$ in $H_{0,rad}^{1}(D)$ with $\lambda > \widetilde{\lambda} > 0$ s.t. $\| u_{\lambda}\|_{L^{2}(D)} = \| u_{\widetilde{\lambda}}\|_{L^{2}(D)} = c$;\\
- (\ref{eq1.1}) admits at least one solution $u_{\lambda}$ in $H_{0,rad}^{1}(D)$ with some $\lambda > 0$ s.t. $\| u_{\lambda}\|_{L^{2}(D)}=\eta_{1}$; \\
- for any $c < \eta_{1}$, (\ref{eq1.1}) admits no solution $u_{\lambda}$ in $H_{0,rad}^{1}(D)$ s.t. $\| u_{\lambda}\|_{L^{2}(D)} = c$.

$(ii)$ If $N = 2$ and $p = 6$, there exists $\eta_{2} > 0$ such that:\\
- For any $c > \eta_{2}$, (\ref{eq1.1}) admits at least one solution $u_{\lambda}$ in $H_{0,rad}^{1}(D)$ with some $\lambda > 0$ s.t. $\| u_{\lambda}\|_{L^{2}(D)} = c$; \\
- for any $c < \eta_{2}$,(\ref{eq1.1}) admits no solution $u_{\lambda}$ in $H_{0,rad}^{1}(D)$ s.t. $\| u_{\lambda}\|_{L^{2}(D)} = c$.

$(iii)$ If $N = 2$ and $p > 6$, for any $c > 0$, (\ref{eq1.1}) admits at least one solution $u_{\lambda}$ in $H_{0,rad}^{1}(D)$ with some $\lambda > 0$ s.t. $\| u_{\lambda}\|_{L^{2}(D)} = c$.

$(iv)$ If $N = 2$ and $p = 4$, there exists $\eta_{3} > 0$ such that:\\
- For any $c > \eta_{3}$, (\ref{eq1.1}) admits at least one solution $u_{\lambda}$ in $H_{0,rad}^{1}(D)$ with some $\lambda > 0$ s.t. $\| u_{\lambda}\|_{L^{2}(D)} = c$; \\
- for any $c < \eta_{3}$,(\ref{eq1.1}) admits no solution $u_{\lambda}$ in $H_{0,rad}^{1}(D)$ s.t. $\| u_{\lambda}\|_{L^{2}(D)} = c$.

$(v)$ If $2 < p < 2 + 4/N$, for any $c > 0$, (\ref{eq1.1}) admits at least one solution $u_{\lambda}$ in $H_{0,rad}^{1}(D)$ with some $\lambda > 0$ s.t. $\| u_{\lambda}\|_{L^{2}(D)} = c$.
\end{theorem}

\begin{remark}
  The thresholds $\eta_1, \eta_2$ and $\eta_3$ are explicitly given in the proof of the Theorem.
\end{remark}

The following table illustrates our main result:

\begin{center}
\begin{tabular}{|m{1.2cm}|m{3.5cm}|m{2.4cm}|m{2.4cm}|m{2.4cm}|m{2.4cm}|m{1.3cm}|}
  \hline
   & $N\geq3$, $p_c \leq p<2^*$ or $N=2,$ $4<p<6$& $N=2,$ $p=6$ & $N=2,$ $p>6$ & $N=2,$ $p=4$ & $2<p<p_c$ \\
   \hline
  $c>\eta_1$ &  at least two normalized solutions &  &  &  &  \\
  \hline
  $c=\eta_1$ & at least one normalized solution &  &  &  &  \\
  \hline
  $c<\eta_1$ & no normalized solution  &  &  &  &  \\
  \hline
  $c>\eta_2$ &  & at least one normalized solution &  &  &  \\
  \hline
  $c<\eta_2$ &  & no normalized solution &  &  &  \\
  \hline
  $c>\eta_3$ &  &  &  & at least one normalized solution  &  \\
  \hline
  $c<\eta_3$ &  &  &  & no normalized solution  &  \\
  \hline
  $c>0$ &  &  & at least one normalized solution &  & at least one normalized solution \\
  \hline
\end{tabular}
\end{center}

\medskip
For any $\lambda > 0$, \cite[Theorem 1.1]{FMT} shows that (\ref{eq1.1}) possesses a unique and non-degenerate solution $u_{\lambda}$ in $H_{0,rad}^{1}(D)$, the Sobolev space of radially symmetric functions. Solutions of (\ref{eq1.1}) can be obtained as critical points of the associated action functional
\begin{equation}
\Phi_{\lambda}(u) := \frac{1}{2}\int_{D} (|\nabla u|^{2} + \lambda u^{2}) - \frac{1}{p}\int_{D}|u|^{p}.
\nonumber
\end{equation}
In fact, $u_{\lambda}$ is the ground state of (\ref{eq1.1}), i.e. $\Phi_{\lambda}(u_{\lambda}) = \inf_{u \in N_{\lambda}}\Phi_{\lambda}(u),$ where $N_{\lambda}$ is the Nehari manifold with a parameter $\lambda$ defined by
\begin{equation}
N_{\lambda}:= \{u \in H_{0,rad}^{1}(D) \backslash \{0\}: \int_{D}| \nabla u |^{2} + \lambda u^{2} - |u|^{p} = 0\}.
\nonumber
\end{equation}
Using \cite[Corollary 2.5]{Song}, we know that $\{(\lambda,u_{\lambda}): \lambda > 0\}$ is a continuous curve in $\mathbb{R} \times H_{0,rad}^{1}(D)$. We will analyze the behavior of $\int_{D}u_{\lambda}^{2}$ when $\lambda \rightarrow +\infty$ and $\lambda \rightarrow 0$ to complete the proof of Theorem \ref{thm1.2}.

As a byproduct of our method, we can obtain the orbital stability/instability of standing wave solutions $e^{i\lambda t}u_{\lambda}(x)$ for (\ref{eq1.3}). Recall that such solutions are called orbitally stable if for each $\epsilon > 0$, there exists $\delta > 0$ such that, whenever $\Psi_{0} \in H_{0,rad}^{1}(D, \mathbb{C})$ is such that $\|\Psi_{0} - u_{\lambda}\|_{H_{0}^{1}(D, \mathbb{C})} < \delta$ and $\Psi(t, x)$ is the solution of (\ref{eq1.3}) with $\Psi(0, \cdot) = \Psi_{0}$ in some interval $[0, t_{0})$, then $\Psi(t, \cdot)$ can be continued to a solution in $0 \leq t < \infty$ and
\[
\sup_{0 < t < \infty}\inf_{s \in \mathbb{R}} \|\Psi(t, x) - e^{i\lambda s}u_{\lambda}(x)\|_{H_{0}^{1}(D, \mathbb{C})} < \epsilon;
\]
otherwise, they are called unstable. We will assume the following condition:

$(LWP)$ For each $\Psi_{0} \in H_{0}^{1}(D, \mathbb{C})$, there exist $t_{0} > 0$, only depending on $\|\Psi_{0}\|_{H_{0}^{1}}$, and a unique solution $\Psi(t, x)$ of (\ref{eq1.3}) with initial datum $\Psi_{0}$ in the interval $I = [0, t_{0})$.

\begin{theorem}
\label{thm1.3} $(i)$ If $N \geq 3,$ and $2 + 4/N \leq p < 2^{\ast}$ or if $N = 2,$ and $4 < p < 6$. Let $u_{\lambda}$ and $u_{\widetilde{\lambda}}$ be given by Theorem \ref{thm1.2} $(i)$, and $\Psi(t, x) = e^{i\lambda t}u_{\lambda}(x)$, $\widetilde{\Psi}(t, x) = e^{i\widetilde{\lambda} t}u_{\widetilde{\lambda}}(x)$. Then for a.e. $c \in (\eta_1, +\infty)$, $\Psi$ is orbitally stable while $\widetilde{\Psi}$ is orbitally unstable in $H_{0,rad}^{1}(D, \mathbb{C})$.

$(ii)$ If $N = 2$ and $p = 6$. Let $u_{\lambda}$ be given by Theorem \ref{thm1.2} $(ii)$, and $\Psi(t, x) = e^{i\lambda t}u_{\lambda}(x)$. Then for a.e. $c \in (\eta_2, +\infty)$, $\Psi$ is orbitally unstable in $H_{0,rad}^{1}(D, \mathbb{C})$.

$(iii)$ If $N = 2$ and $p > 6$. Let $u_{\lambda}$ be given by Theorem \ref{thm1.2} $(iii)$, and $\Psi(t, x) = e^{i\lambda t}u_{\lambda}(x)$. Then for a.e. $c \in (0, +\infty)$, $\Psi$ is orbitally unstable in $H_{0,rad}^{1}(D, \mathbb{C})$.

$(iv)$ If $N = 2$ and $p = 4$. Let $u_{\lambda}$ be given by Theorem \ref{thm1.2} $(iv)$, and $\Psi(t, x) = e^{i\lambda t}u_{\lambda}(x)$. Then for a.e. $c \in (\eta_3, +\infty)$, $\Psi$ is orbitally stable in $H_{0,rad}^{1}(D, \mathbb{C})$.

$(v)$ If $2 < p < 2 + 4/N$. Let $u_{\lambda}$ be given by Theorem \ref{thm1.2} $(v)$, and $\Psi(t, x) = e^{i\lambda t}u_{\lambda}(x)$. Then for a.e. $c \in (0, +\infty)$, $\Psi$ is orbitally stable in $H_{0,rad}^{1}(D, \mathbb{C})$.
\end{theorem}

It is worth noting that something interesting happens when $p \geq 2 + 4/N$. If we consider (\ref{eq1.3}) and (\ref{eq1.1}) on the entire space $\mathbb{R}^{N}$, for each $c > 0$, (\ref{eq1.1}) has a unique solution with prescribed mass $c$ and the corresponding standing wave is orbitally unstable. The nature of the domain can affect the number of normalized solutions and have a stabilizing effect. More precisely, the same problem on the unit ball when $p > 2 + 4/N$, admits at least two normalized solutions. The one with a small mass is orbitally stable. For large masses $c$, there exists no normalized solution. Such results can be found in \cite{FM,FSK,NTV,Song3,Song2}. We also establish a new critical exponent $p = 6$ when $N = 2$.  If $N \geq 3,$ $2 + 4/N \leq p < 2^{\ast}$ or if $N = 2,$ $4 < p < 6$, contrary to the case of the unit ball case, there exist no normalized solution with a small mass, while there exist at least two normalized solutions with a large mass such as the one with a larger $\lambda$ is orbitally stable. For problems on the unit ball or $\mathbb{R}^{N}$, radial solutions of (\ref{eq1.1}) have a unique maximum point while maximum points of radial solutions for (\ref{eq1.1}) on $D$ are a circle. The exponent $p = 6$ comes from this difference. For this reason, we believe that if $N \geq 3,$ and $2 < p < 2^{\ast}$  or if $N = 2,$ and $2 < p < 6$, (\ref{eq1.1}) in an annulus possesses a normalized solution for each $c > 0$ and the corresponding standing wave is orbitally stable (which will be addressed in another work).

Let $u_{\lambda}$ be the unique solution of (\ref{eq1.1}) in $H_{0,rad}^{1}(D)$. Then
$$
u_{\lambda,R}(x) = R^{\frac{2}{2-p}}u_{\lambda R^{2}}(\frac{x}{R})
$$
is the unique solution in $H_{0,rad}^{1}(B_{R}^{c})$ of
\begin{equation} \label{eq1.4}
\left\{
\begin{array}{cc}
-\Delta u + \lambda u = u^{p-1}, x \in B_{R}^{c}, \\
u > 0, u = 0 \ on \ \partial B_{R}, u(x) \rightarrow 0 \ as \ |x| \rightarrow +\infty,
\end{array}
\right.
\end{equation}
As a corollary of Theorem \ref{thm1.2}, we have:

\begin{corollary}
\label{cor1.4} $(i)$ If $N \geq 3,$ and $2 + 4/N \leq p < 2^{\ast}$ or if $N = 2$ and $4 < p < 6$, there exists some $\eta_{1,R} > 0$ such that for any $c > \eta_{1,R}$, (\ref{eq1.4}) admits at least two solutions $u_{\lambda}$ and $u_{\widetilde{\lambda}}$ with $\lambda > \widetilde{\lambda} > 0$ s.t. $\| u_{\lambda}\|_{L^{2}(B_{R}^{c})} = \| u_{\widetilde{\lambda}}\|_{L^{2}(B_{R}^{c})} = c$; (\ref{eq1.4}) admits at least one solution $u_{\lambda}$ with some $\lambda > 0$ s.t. $\| u_{\lambda}\|_{L^{2}(B_{R}^{c})} = \eta_{1,R}$; for any $c < \eta_{1,R}$, (\ref{eq1.4}) admits no solution $u_{\lambda}$ s.t. $\| u_{\lambda}\|_{L^{2}(B_{R}^{c})} = c$.

$(ii)$ If $N = 2$ and $p = 6$, there exists some $\eta_{2,R} > 0$ such that for any $c > \eta_{2,R}$, (\ref{eq1.4}) admits at least one solution $u_{\lambda}$ with some $\lambda > 0$ s.t. $\| u_{\lambda}\|_{L^{2}(B_{R}^{c})} = c$; for any $c < \eta_{2,R}$, (\ref{eq1.4}) admits no solution $u_{\lambda}$ s.t. $\| u_{\lambda}\|_{L^{2}(B_{R}^{c})} = c$.

$(iii)$ If $N = 2$ and $p > 6$, for any $c > 0$, (\ref{eq1.4}) admits at least one solution $u_{\lambda}$ with some $\lambda > 0$ s.t. $\| u_{\lambda}\|_{L^{2}(B_{R}^{c})} = c$.

$(iv)$ If $N = 2$ and $p = 4$, there exists some $\eta_{3,R} > 0$ such that for any $c > \eta_{3,R}$, (\ref{eq1.4}) admits at least one solution $u_{\lambda}$ with some $\lambda > 0$ s.t. $\| u_{\lambda}\|_{L^{2}(B_{R}^{c})} = c$; for any $c < \eta_{3,R}$, (\ref{eq1.4}) admits no solution $u_{\lambda}$ s.t. $\| u_{\lambda}\|_{L^{2}(B_{R}^{c})} = c$.

$(v)$ If $2 < p < 2 + 4/N$, for any $c > 0$, (\ref{eq1.4}) admits at least one solution $u_{\lambda}$ with some $\lambda > 0$ s.t. $\| u_{\lambda}\|_{L^{2}(B_{R}^{c})} = c$.
\end{corollary}

For the asymptotic behaviors of $\eta_{1,R}$, $\eta_{2,R}$ and $\eta_{3,R}$, we have the following result.

\begin{theorem}
\label{thm1.5} Let $\eta_{1,R}$, $\eta_{2,R}$ and $\eta_{3,R}$ be given by Corollary \ref{cor1.4}.

$(i)$  If $N \geq 3$ and $2 + 4/N < p < 2^{\ast}$ or if $N = 2$ and $4 < p < 6$,
$$
\lim_{R \rightarrow 0}\eta_{1,R} = 0, \lim_{R \rightarrow +\infty}\eta_{1,R} = +\infty.
$$
In particular, for any given $c > 0$, there exists $R_{1} = R_{1}(c)$ small and $\widetilde{R}_{1} = \widetilde{R}_{1}(c)$ large, such that for any $R < R_{1}$, (\ref{eq1.4}) admits at least two solutions with prescribed mass $c$ in $H_{0,rad}^{1}(B_{R}^{c})$ while for any $R > \widetilde{R}_{1}$, (\ref{eq1.4}) admits no solution in $H_{0,rad}^{1}(B_{R}^{c})$ with prescribed mass $c$.

$(ii)$ If $N = 2$ and $p = 6$,
$$
\lim_{R \rightarrow 0}\eta_{2,R} = 0, \lim_{R \rightarrow +\infty}\eta_{2,R} = +\infty.
$$
In particular, for any given $c > 0$, there exists $R_{2} = R_{2}(c)$ small and $\widetilde{R}_{2} = \widetilde{R}_{2}(c)$ large, such that for any $R < R_{2}$, (\ref{eq1.4}) admits at least one solution in $H_{0,rad}^{1}(B_{R}^{c})$ with prescribed mass $c$ while for any $R > \widetilde{R}_{2}$, (\ref{eq1.4}) admits no solution in $H_{0,rad}^{1}(B_{R}^{c})$ with prescribed mass $c$.

$(iii)$ If $N \geq 3$ and $p = 2 + 4/N$, then $\eta_{1,R} \equiv \eta_{1}$; when $N = 2$ and $p = 4$, then $\eta_{3,R} \equiv \eta_{3}$.
\end{theorem}

We organize this paper as follows. In Section 2, we study the asymptotic behaviors when $\lambda \rightarrow +\infty$. A new Pohozaev-type identity without boundary will be established in this section. In Section 3, we analyze the asymptotic behaviors when $\lambda \rightarrow 0$ and complete the proofs of Theorems \ref{thm1.2} and \ref{thm1.5}. Finally in Section 4, we discuss the orbital stability/instability of standing wave solutions and show Theorem \ref{thm1.3}.

\section{Asymptotic behaviors when $\lambda \rightarrow +\infty$}

\begin{lemma}
\label{lemA.1} Let $u_{\lambda} = u_{\lambda}(r)$ be the unique solution of (\ref{eq1.1}) in $H_{0,rad}^{1}(D)$. Then $u_{\lambda}(r)$ has a unique maximum point $\bar{r}_{\lambda}$ and $u_{\lambda}' > 0$ in $(1,\bar{r}_{\lambda})$, $u_{\lambda}' < 0$ in $(\bar{r}_{\lambda},+\infty)$.
\end{lemma}

\textit{Proof.  } Standard arguments show that $u_{\lambda}(r) \in C^{2}([1,+\infty))$. By the strong maximum principle, $u_{\lambda}'(1) > 0$. Since $u_{\lambda}(r) \rightarrow 0$ as $r \rightarrow +\infty$, we can assume that $u_{\lambda}'(\bar{r}_{\lambda}) = 0$ and $u_{\lambda}' > 0$ in $(1,\bar{r}_{\lambda})$. Then $u_{\lambda}$ satisfies the following equation
\begin{equation} \label{eqa.1}
\left\{
\begin{array}{cc}
-(u'' + \frac{N-1}{r}u') + \lambda u = u^{p-1} \ in \ (\bar{r}_{\lambda},+\infty), \\ u(\bar{r}_{\lambda}) = \alpha, u'(\bar{r}_{\lambda}) = 0, u(r) \rightarrow 0 \ as \ r \rightarrow +\infty,
\end{array}
\right.
\end{equation}
where $\alpha = u_{\lambda}(\bar{r}_{\lambda})$. By \cite[Theorem 2]{Sir}, $u'(r) < 0$ in $(\bar{r}_{\lambda},+\infty)$. \qed\vskip 5pt

According to Lemma \ref{lemA.1}, if $\bar{r}_{\lambda}$ be the unique maximum point of $u_{\lambda},$ we have two possibilities:

$(A_{1})$ $\liminf_{\lambda \rightarrow +\infty}\bar{r}_{\lambda} < +\infty$.

$(A_{2})$ $\bar{r}_{\lambda} \rightarrow +\infty$ as $\lambda \rightarrow +\infty$.

Our first main result in this section reads as follows:

\begin{theorem}
\label{thm2.2} Let $u_{\lambda}$ be the unique solution of (\ref{eq1.1}) in $H_{0,rad}^{1}(D)$ with $\lambda \rightarrow +\infty$.

$(i)$ When $N = 2$ and $p < 6$ or $N \geq 3$, $\int_{D}u_{\lambda}^{2} \rightarrow +\infty$.

$(ii)$ When $N = 2$ and $p = 6$, $\liminf_{\lambda \rightarrow +\infty}\int_{D}u_{\lambda}^{2} > 0$. If we also assume that $(A_{1})$ holds, then $\liminf_{\lambda \rightarrow +\infty}\int_{D}u_{\lambda}^{2} < +\infty$.

$(iii)$ Assume that $(A_{2})$ holds. When $N = 2$ and $p = 6$, $\int_{D}u_{\lambda}^{2} \rightarrow +\infty$.

$(iv)$ Assume that $(A_{1})$ holds. When $N = 2$ and $p > 6$, $\liminf_{\lambda \rightarrow +\infty}\int_{D}u_{\lambda}^{2} = 0$.
\end{theorem}

\begin{lemma}
\label{lemA.2} Let $u_{\lambda}$ be the unique solution of (\ref{eq1.1}) in $H_{0,rad}^{1}(D)$. Then
\begin{equation}
0 < \liminf_{\lambda \rightarrow +\infty}\frac{\lambda}{\|u_{\lambda}\|_{L^{\infty}}^{p-2}} \leq \limsup_{\lambda \rightarrow +\infty}\frac{\lambda}{\|u_{\lambda}\|_{L^{\infty}}^{p-2}} \leq 1.
\end{equation}
\end{lemma}

\textit{Proof.  } Let $\bar{r}_{\lambda}$ be the unique maximum point of $u_{\lambda}$. Then
\begin{equation}
\lambda u_{\lambda}(\bar{r}_{\lambda}) < -u_{\lambda}''(\bar{r}_{\lambda}) + \lambda u_{\lambda}(\bar{r}_{\lambda}) = u_{\lambda}(\bar{r}_{\lambda})^{p-1},
\end{equation}
i.e. $\lambda < u_{\lambda}(\bar{r}_{\lambda})^{p-2}$.

Next we will show that $\liminf_{\lambda \rightarrow +\infty}\lambda/\|u_{\lambda}\|_{L^{\infty}}^{p-2} > 0$. We will argue by contradiction: Suppose on the contrary that there exists $\{\lambda_{n}\}$ with $\lambda_{n} \rightarrow +\infty$ such that $\lambda_{n}/\|u_{\lambda_{n}}\|_{L^{\infty}}^{p-2} \rightarrow 0$. To simplify the notation, we shall denote $u_{n} = u_{\lambda_{n}}$ and $\bar{r}_{n} = \bar{r}_{\lambda_{n}}$. Notice that $u_{n}$ satisfies
\begin{equation}
-(u_{n}'' + \frac{N-1}{r}u_{n}') + \lambda_{n} u_{n} = u_{n}^{p-1} \ in \ I, u_{n}(1) = 0,
\end{equation}
where $I = (1,+\infty)$. We consider
\begin{equation} \label{eqa.7}
v_{n}(r) = \frac{1}{\|u_{n}\|_{L^{\infty}(I)}} u_{n}(\frac{r}{\|u_{n}\|_{L^{\infty}(I)}^{\frac{p-2}{2}}} + \bar{r}_{n}), r \in I_{n} = (\|u_{n}\|_{L^{\infty}(I)}^{\frac{p-2}{2}}(1 - \bar{r}_{n}),+\infty).
\end{equation}
Then $\|v_{n}\|_{L^{\infty}(I_{n})} = v_{n}(0) = 1$ and
\begin{equation} \label{eqa.3}
-(v_{n}'' + \frac{N-1}{r + \|u_{n}\|_{L^{\infty}(I)}^{\frac{p-2}{2}}\bar{r}_{n}}v_{n}') + \frac{\lambda_{n}}{\|u_{n}\|_{L^{\infty}}^{p-2}}v_{n} = v_{n}^{p-1} \ in \ I_{n}, v_{n} = 0 \ on \ \partial I_{n}.
\end{equation}

Assume that $\bar{r}_{n} \rightarrow \bar{r} \in \overline{I} \cup \{+\infty\}$ up to a subsequence. ($\bar{r} = +\infty$ means that $\bar{r}_{n} \rightarrow +\infty$.) If
$$
\lim_{n \rightarrow +\infty}\|u_{n}\|_{L^{\infty}(I)}^{\frac{p-2}{2}}(1 - \bar{r}_{n}) = -\infty,
$$
$I_{n} \rightarrow \mathbb{R}$. For any $R > 0$, $(-R,R)$ is contained in $I_{n}$ for $n$ large enough. Applying standard elliptic estimates and passing to a subsequence if necessary, we may assume that $v_{n} \rightarrow v$ in $C_{loc}(\mathbb{R})$ where $v(0) = 1$ and $v$ satisfies
\begin{equation} \label{eqa.4}
-v'' = v^{p-1} \ in \ \mathbb{R}, 0 \leq v \leq 1 \ in \ \mathbb{R}, v \in C^{2}(\mathbb{R}).
\end{equation}
If
$$
\lim_{n \rightarrow +\infty}\|u_{n}\|_{L^{\infty}(I)}^{\frac{p-2}{2}}(1 - \bar{r}_{n}) = -k \leq 0,
$$
$I_{n} \rightarrow (-k,+\infty)$. Similarly, we obtain a solution $v$ with $v(0) = 1$ of
\begin{equation} \label{eqa.5}
-v'' = v^{p-1} \ in \ (-k,+\infty), v(-k) = 0, 0 \leq v \leq 1, v \in C^{2}([-k,+\infty)).
\end{equation}
Noticing that $v(0) = 1$ and $v(-k) = 0$, we know that $k > 0$. Then the following Lemma \ref{lemA.3} provides a contradiction. \qed\vskip 5pt

\begin{lemma}
\label{lemA.3} (\ref{eqa.4}) (resp. (\ref{eqa.5})) admits no solution $v$ that satisfies $v(0) = 1$.
\end{lemma}

\textit{Proof.  } If $v$ is such that $v(0) = 1$ is a solution of (\ref{eqa.4}) (resp. (\ref{eqa.5})), then $v$ satisfies the following equation
\begin{equation} \label{eqa.6}
-v'' = v^{p-1}, v(0) = 1, v'(0) = 0.
\end{equation}
Let
$$
E(r) = \frac{1}{2}(v'(r))^{2} + \frac{1}{p}(v(r))^{p}.
$$
Since
$$
E'(r) = v'(r)v''(r) + (v(r))^{p-1}v'(r) = 0,
$$
$E(r) \equiv E(0) = 1/p$, implying that $v(r) \leq 1$. $-v''(0) = 1 > 0$, which implies that $v$ cannot be a constant and $v'(r) < 0$ for small positive $r$.

We claim that there exists $\widehat{r} > 0$ such that $v(\widehat{r}) = 0$. If this is not the case, $0 < v(r) \leq 1$ for any $r > 0$. If there exists $\widetilde{r} > 0$ such that $v'(\widetilde{r}) = 0$, we may assume that $v'(r) < 0$ in $(0,\widetilde{r})$. Then $v(\widetilde{r}) < v(0) = 1$, implying that
$$
E(\widetilde{r}) = \frac{1}{p}(v(\widetilde{r}))^{p} < \frac{1}{p},
$$
which is a contradiction with $E(r) \equiv 1/p$. Thus $v'(r) < 0$ for any $r > 0$. Since $v > 0$ in $(0,+\infty)$, we have $\lim_{r \rightarrow +\infty}v'(r) = 0$. We assume that $\lim_{r \rightarrow +\infty}v(r) = a \in [0,1)$. Then
$$
\lim_{r \rightarrow +\infty}E(r) = \frac{1}{p}a^{p} < \frac{1}{p},
$$
contradicting the fact that $E(r) \equiv 1/p$.

Therefore, the proof of the claim is now complete and we may assume that $v(\widehat{r}) = 0,$ and $v(r) > 0$ in $(0,\widehat{r})$. Furthermore, $v'(\widehat{r}) \neq 0$, implying that $v'(\widehat{r}) < 0$. Thus $v(\widehat{r} + \epsilon) < 0$ for $\epsilon > 0$ small enough, which contradicts the fact that $v \geq 0$. \qed\vskip 5pt

\begin{theorem}
\label{thmA.4} Let $u_{n}(r) = u_{\lambda_{n}}$ be the unique solution of (\ref{eq1.1}) in $H_{0,rad}^{1}(D)$ with $\lambda = \lambda_{n} \rightarrow +\infty$ and
\begin{equation}
\omega_{n}(r) = \lambda_{n}^{\frac{1}{2-p}} u_{n}(\frac{r}{\sqrt{\lambda_{n}}} + \bar{r}_{n}), r \in \widetilde{I}_{n} = (\sqrt{\lambda_{n}}(1 - \bar{r}_{n}),+\infty),
\end{equation}
where $\bar{r}_{n} = \bar{r}_{\lambda_{n}}$ is the unique maximum point of $u_{n}$. Passing to a subsequence if necessary, we may assume that $\omega_{n}(r) \rightarrow W(r)$ in $H_{0}^{1}(\mathbb{R})$ where $W \in H^{1}(\mathbb{R})$ is the unique solution of
\begin{equation} \label{eqa.10}
-W'' + W = W^{p-1} \ in \ \mathbb{R}, W(r) \rightarrow 0 \ as \ |r| \rightarrow +\infty.
\end{equation}
\end{theorem}

\textit{Proof.  } Note that $\omega_{n}$ satisfies
\begin{equation}
-(\omega_{n}'' + \frac{N-1}{r + \sqrt{\lambda_{n}}\bar{r}_{n}}\omega_{n}') + \omega_{n} = \omega_{n}^{p-1} \ in \ \widetilde{I}_{n}, \omega_{n} = 0 \ on \ \partial \widetilde{I}_{n}.
\end{equation}
By Lemma \ref{lemA.2}, $\omega_{n}(0) = \|\omega_{n}\|_{L^{\infty}}$ is bounded. We claim that
$$
\lim_{n \rightarrow +\infty}\sqrt{\lambda_{n}}(1 - \bar{r}_{n}) = -\infty.
$$
If not, up to a subsequence, we may assume that
$$
\lim_{n \rightarrow +\infty}\sqrt{\lambda_{n}}(1 - \bar{r}_{n}) = -k \leq 0,
$$
Similar to the proof of Lemma \ref{lemA.2}, we obtain a solution $W$ with $W(0) = \max W(r)$ of
\begin{equation}
-W'' + W = W^{p-1} \ in \ (-k,+\infty), W(-k) = 0, W \geq 0, W \in C^{2}([-k,+\infty)).
\end{equation}
Notice that $W(k) = W(-k) = 0$ and $W'(k) = -W'(-k) < 0$, thus we find a contradiction with $W \geq 0$. Therefore, the claim has been proved, and we know that $\omega_{n} \rightarrow W$ in $C_{loc}^{2}(\mathbb{R})$ where $W(0) = \max W(r)$ and $W$ satisfies
\begin{equation}
-W'' + W = W^{p-1} \ in \ \mathbb{R}, W \geq 0, W \in C^{2}(\mathbb{R}).
\end{equation}
According to Lemma \ref{lemA.1}, $W$ is increasing in $(-\infty,0)$ and decreasing in $(0,+\infty)$. Thus $W(r) \rightarrow 0$ as $r \rightarrow +\infty$, i.e. $W$ is the unique solution of (\ref{eqa.10}).

Next we will show that $\omega_{n}$ decays uniformly and $\omega_{n} \rightarrow W$ in $H_{0}^{1}(\mathbb{R})$. Suppose on the contrary that there exists $r_{n}$ with $|r_{n}| \rightarrow +\infty$ such that $\omega_{n}(r_{n}) = \delta$ where $\delta > 0$ can be chosen small enough. Set $\widetilde{\omega}_{n}(r) = \omega_{n}(r + r_{n})$.

If $r_{n} \rightarrow +\infty$, then $\widetilde{\omega}_{n}$ converges to $\widetilde{W}$, a nontrivial solution of the
following equation
\begin{equation} \label{eqa.8}
-\widetilde{W}'' + \widetilde{W} = \widetilde{W}^{p-1} \ in \ \mathbb{R},
\end{equation}
with $\widetilde{W} \geq 0$ and bounded, $\widetilde{W}(0) = \delta$. According to Lemma \ref{lemA.1}, $\widetilde{\omega}_{n}$ is decreasing in $[-r_{n},+\infty)$ and thus $\widetilde{W}$ is decreasing in $\mathbb{R}$. By Lemma \ref{lemA.5} (see below), (\ref{eqa.8}) admits a no nontrivial decreasing solution, which is a contradiction.

If $r_{n} \rightarrow -\infty$ and $\sqrt{\lambda_{n}}(1 - \bar{r}_{n}) - r_{n} \rightarrow -\infty$, we can find a contradiction by using similar arguments as above. If $\sqrt{\lambda_{n}}(1 - \bar{r}_{n}) - r_{n} \nrightarrow -\infty$, passing to a subsequence if necessary, we may assume that $\sqrt{\lambda_{n}}(1 - \bar{r}_{n}) - r_{n} \rightarrow -k \leq 0$. Then $\widetilde{\omega}_{n}$ converges to $\widetilde{W}$, a nontrivial solution of the
following equation
\begin{equation} \label{eqa.9}
-\widetilde{W}'' + \widetilde{W} = \widetilde{W}^{p-1} \ in \ (-k,+\infty), \widetilde{W}(-k) = 0,
\end{equation}
with $\widetilde{W} \geq 0$ and bounded, $\widetilde{W}(0) = \delta$. Noticing that $\widetilde{W}(-k) = 0$, we know that $k > 0$. According to Lemma \ref{lemA.1}, $\widetilde{\omega}_{n}$ is increasing in $(\sqrt{\lambda_{n}}(1 - \bar{r}_{n}) - r_{n},-r_{n}]$ and thus $\widetilde{W}$ is increasing in $(-k,+\infty)$. By the following Lemma \ref{lemA.5}, (\ref{eqa.9}) admits a no nontrivial increasing solution, which contradicts our assumption.

We now show the uniform decay for $\omega_{n}$, i.e. for any $\epsilon > 0$, there exists $R > 0$ such that for all $r$ with $|r| > R$, $\omega_{n}(r) < \epsilon$. Notice that
$$
r + \sqrt{\lambda_{n}}\bar{r}_{n} > \sqrt{\lambda_{n}}(1 - \bar{r}_{n}) + \sqrt{\lambda_{n}}\bar{r}_{n} = \sqrt{\lambda_{n}} \rightarrow +\infty.
$$
Take $\epsilon$ small enough, we have
\begin{equation}
-(\omega_{n}'' + \frac{N-1}{r + \sqrt{\lambda_{n}}\bar{r}_{n}}\omega_{n}') + \frac{1}{2}\omega_{n} \leq 0, \forall |r| > R.
\end{equation}
Then there exist $C > 0$ and $K > 0$ such that for all $n > K$, $|r| > R$
$$
\omega_{n}(r) \leq C e^{-\frac{r}{2}}.
$$
Combining the convergence of $\omega_{n}$ in $C_{loc}^{2}(\mathbb{R})$ and the uniform exponential decay, we know that $\omega_{n}(r) \rightarrow W(r)$ in $H_{0}^{1}(\mathbb{R})$ where $W \in H_{0}^{1}(\mathbb{R})$ is the unique solution of (\ref{eqa.10}).
\qed\vskip 5pt

\begin{lemma}
\label{lemA.5} Let $W \in C^{2}$ satisfying $-W'' + W = |W|^{p-2}W$ in $\mathbb{R}$. Unless $W \equiv 0$ or $W \equiv\pm 1$, $W$ is not monotone. Furthermore, if $W$ has a zero point, then $W$ is periodic.
\end{lemma}

\textit{Proof.  } Step 1: $W \in L^{\infty}(\mathbb{R})$.

Let
$$
E(r) = \frac{1}{2}(W'(r))^{2} + \frac{1}{p}|W|^{p} - \frac{1}{2}W^{2}.
$$
Then $E'(r) = 0$ and thus $E(r)$ is a constant. Therefore,
$$
\limsup_{r \rightarrow +\infty}(\frac{1}{p}|W|^{p} - \frac{1}{2}W^{2}) \leq \lim_{r \rightarrow +\infty}E(r) = E(0) < +\infty,
$$
implying the boundedness of $|W(r)|$ when $r \rightarrow +\infty$.

Step 2: If $W$ possesses an extreme point, then $W \equiv 0$ or $W \equiv \pm 1$ or $W$ is not monotone. Furthermore, if $W$ also possesses a zero point, then $W$ is periodic.

Without loss of generality, we may assume that $\max W(r) = W(0) := \alpha$. (If $W$ possesses a minimum point, we consider $-W$.) Then $W$ satisfies
\begin{equation} \label{eqa.21}
-W'' + W = |W|^{p-2}W, W(0) = \alpha, W'(0) = 0.
\end{equation}
If $W(r) \equiv \alpha$, then $\alpha = 0$ or $\alpha = \pm 1$. If $\alpha \neq 0$ and $\alpha \neq \pm 1$, then $W''(0) = W(1 - |W|^{p-2}) \neq 0$. Since $\max W(r) = W(0)$, we know that $W''(0) < 0$ and $W$ is not monotone.

Without loss of generality, we assume that $\alpha > 0$. If $W$ also possesses a zero point, we can assume that $W(\widehat{r}) = 0$ and $W(r) > 0$ in $(0,\widehat{r})$. Let
\begin{equation}
\widetilde{W}(r) =
\left\{
\begin{array}{cc}
\ \ \ \ W(r), \ \ \ \ \ \ \ 0 \leq r \leq \widehat{r}, \\
-W(2\widehat{r} - r), \ \ \widehat{r} < r \leq 2\widehat{r}, \\
\ \ W(4\widehat{r} - r), \ \ 2\widehat{r} < r \leq 4\widehat{r}, \\
\widetilde{W}(r + 4\widehat{r}), \ \ \ r \in \mathbb{R}.
\end{array}
\right.
\end{equation}
Then $\widetilde{W}$ solves (\ref{eqa.21}). By the uniqueness of the solution of (\ref{eqa.21}), $W \equiv \widetilde{W}$ is periodic.

Step 3: $W$ possesses an extreme point.

Otherwise, $W$ is monotone. Without loss of generality, we may assume that $W$ is decreasing on $\mathbb{R}$. By Step 1, we may assume that $\lim_{r \rightarrow -\infty}W(r) = a$. $\lim_{r \rightarrow -\infty}W''(r) = (1 - |a|^{p-2})a < 0$. Let $u$ be the solution of
\begin{equation}
-u'' + u = |u|^{p-2}u, u(0) = a, u'(0) = 0.
\end{equation}
$u''(r) = (1 - |a|^{p-2})a < 0$. Take $\epsilon > 0$ small, there exist $r_{1} > 0$, $r_{2} \in \mathbb{R}$ such that $u(r_{1}) = a-\epsilon$, $u'(r_{1}) < 0$ and $W(r_{2}) = a-\epsilon$. Then
$$
u'(r_{1}) = -\sqrt{\frac{2}{p}|a|^{p} - a^{2} - (\frac{2}{p}|a-\epsilon|^{p} - |a-\epsilon|^{2})} = W'(r_{2}).
$$
Therefore, both $u(r + r_{1})$ and $W(r + r_{2})$ satisfy
\begin{equation}
-v'' + v = |v|^{p-2}v, v(0) = a-\epsilon, v'(0) = -\sqrt{\frac{2}{p}|a|^{p} - a^{2} - (\frac{2}{p}|a-\epsilon|^{p} - |a-\epsilon|^{2})},
\end{equation}
implying $u(r + r_{1}) \equiv W(r + r_{2})$. By the proof of step 2, $u$ is not monotone, contradicting the fact that $W$ is monotone. \qed\vskip 5pt

\begin{corollary}
$\omega_{n}$ and $W$ are given by Theorem \ref{thmA.4}. Passing to a subsequence if necessary, we have
\begin{equation} \label{eqa.11}
\int_{\sqrt{\lambda_{n}}(1 - \bar{r}_{n})}^{+\infty}\omega_{n}^{2}r^{k} \rightarrow \int_{-\infty}^{+\infty}W^{2}r^{k}, k = 0, 1, \cdots, N-1.
\end{equation}
\end{corollary}

\textit{Proof.  } (\ref{eqa.11}) directly follows by the convergence of $\omega_{n}$ in $C_{loc}^{2}(\mathbb{R})$ and the uniform exponential decay.
\qed\vskip 5pt

\textbf{Proof of Theorem \ref{thm2.2}:  } Let $u_{n} = u_{\lambda_{n}}$ be the unique solution of (\ref{eq1.1}) in $H_{0,rad}^{1}(D)$ with $\lambda = \lambda_{n} \rightarrow +\infty$ and $\bar{r}_{n} = \bar{r}_{\lambda_{n}}$ be the unique maximum point of $u_{n}$. Direct computations imply that:
\begin{eqnarray} \label{eqa.12}
\int_{D}u_{n}^{2} &=& C\int_{1}^{+\infty}u_{n}^{2}r^{N-1}dr \nonumber \\
&=& C\int_{1 - \bar{r}_{n}}^{+\infty}(u_{n}(r+\bar{r}_{n}))^{2}(r+\bar{r}_{n})^{N-1}dr \nonumber \\
&=& C\int_{1 - \bar{r}_{n}}^{+\infty}(u_{n}(r+\bar{r}_{n}))^{2}(\Sigma_{k = 0}^{N-1}C_{N-1}^kr^k\bar{r}_{n}^{N-1-k})dr \nonumber \\
&=& C\Sigma_{k = 1}^{N-1}C_{N-1}^k\lambda_{n}^{\frac{2}{p-2} - \frac{k+1}{2}}\bar{r}_{n}^{N-1-k}\int_{\sqrt{\lambda_{n}}(1 - \bar{r}_{n})}^{+\infty}\omega_{n}^{2}r^{k}dr,
\end{eqnarray}
where $C = \int_{S^{N-1}}d\sigma> 0$, $C_{N-1}^k = \frac{(N-1)(N-2)\cdots(N-k)}{k(k-1)\cdots 1}$ if $k \geq 1$ and $C_{N-1}^0 = 1$.

Proof of Theorem \ref{thm2.2} $(i)$: We will argue by contradiction,  suppose that there exists $\{u_{n}\}$ such that $u_{n} = u_{\lambda_{n}}$ is the unique solution of (\ref{eq1.1}) in $H_{0,rad}^{1}(D)$ with $\lambda = \lambda_{n} \rightarrow +\infty$ and
\begin{equation} \label{eqa.13}
\lim_{n \rightarrow +\infty}\int_{D}u_{n}^{2} < +\infty.
\end{equation}
When $N = 2$, we assume that $p < 6$. When $N \geq 3$, $p < 2N/(N-2) \leq 6$. Up to a subsequence, (\ref{eqa.12}) shows that for large $n$,
\begin{equation}
\int_{D}u_{n}^{2} \geq \frac{C}{2}\lambda_{n}^{\frac{2}{p-2} - \frac{1}{2}}\bar{r}_{n}^{N-1}\int_{\sqrt{\lambda_{n}}(1 - \bar{r}_{n})}^{+\infty}\omega_{n}^{2}dr \geq \frac{C}{2}\lambda_{n}^{\frac{2}{p-2} - \frac{1}{2}}\int_{\sqrt{\lambda_{n}}(1 - \bar{r}_{n})}^{+\infty}\omega_{n}^{2}dr \rightarrow +\infty,
\end{equation}
which is a contradiction with (\ref{eqa.13}).

Proof of Theorem \ref{thm2.2} $(ii)$: First, we prove that $\liminf_{\lambda \rightarrow +\infty}\int_{D}u_{\lambda}^{2} > 0$. We will argue by contradiction, and assume that there exists $\{u_{n}\}$ such that $u_{n} = u_{\lambda_{n}}$ is the unique solution of (\ref{eq1.1}) in $H_{0,rad}^{1}(D)$ with $\lambda = \lambda_{n} \rightarrow +\infty$ and
\begin{equation} \label{eqa.14}
\int_{D}u_{n}^{2} \rightarrow 0.
\end{equation}
Up to a subsequence, (\ref{eqa.12}) shows that for large $n$
\begin{equation}
\int_{D}u_{n}^{2} \geq \frac{C}{2}\bar{r}_{n}^{N-1}\int_{\sqrt{\lambda_{n}}(1 - \bar{r}_{n})}^{+\infty}\omega_{n}^{2}dr \geq \frac{C}{2}\int_{\sqrt{\lambda_{n}}(1 - \bar{r}_{n})}^{+\infty}\omega_{n}^{2}dr \rightarrow \frac{C}{2}\int_{-\infty}^{+\infty}W^{2}dr > 0,
\end{equation}
where $W$ is given by Theorem \ref{thmA.4}. Thus we find a contradiction with (\ref{eqa.14}).

Next, we prove that $\liminf_{\lambda \rightarrow +\infty}\int_{D}u_{\lambda}^{2} < +\infty$ under the hypothesis $(A_{1})$. By $(A_{1})$, we may assume that there exists $u_{n} = u_{\lambda_{n}}$ with $\bar{r}_{n} \rightarrow \bar{r} \in [1,+\infty)$ where $\bar{r}_{n}$ is the unique maximum point of $u_{n}$.
Up to a subsequence, (\ref{eqa.12}) shows that
\begin{equation}
\int_{D}u_{n}^{2} = C\Sigma_{k = 0}^{N-1}C_{N-1}^k\lambda_{n}^{\frac{1}{2} - \frac{k+1}{2}}\bar{r}_{n}^{N-1-k}\int_{\sqrt{\lambda_{n}}(1 - \bar{r}_{n})}^{+\infty}\omega_{n}^{2}r^{k}dr \rightarrow C\bar{r}^{N-1}\int_{-\infty}^{+\infty}W^{2}dr < +\infty
\end{equation}
where $W$ is given by Theorem \ref{thmA.4}. Consequently we arrive at the conclusion as desired.

Proof of Theorem \ref{thm2.2} $(iii)$:  We will argue by contradiction, let us assume that there exists $\{u_{n}\}$ such that $u_{n} = u_{\lambda_{n}}$ is the unique solution of (\ref{eq1.1}) in $H_{0,rad}^{1}(D)$ with $\lambda = \lambda_{n} \rightarrow +\infty$ and
\begin{equation} \label{eqa.15}
\lim_{n \rightarrow +\infty}\int_{D}u_{n}^{2} < +\infty.
\end{equation}
By $(A_{2})$, $\bar{r}_{n} \rightarrow +\infty$ where $\bar{r}_{n}$ is the unique maximum point of $u_{n}$. Up to a subsequence, (\ref{eqa.12}) shows that for large $n$,
\begin{equation}
\int_{D}u_{n}^{2} = C\Sigma_{k = 0}^{N-1}C_{N-1}^k\lambda_{n}^{\frac{1}{2} - \frac{k+1}{2}}\bar{r}_{n}^{N-1-k}\int_{\sqrt{\lambda_{n}}(1 - \bar{r}_{n})}^{+\infty}\omega_{n}^{2}r^{k}dr \geq \frac{C}{2}\bar{r}_{n}^{N-1}\int_{\sqrt{\lambda_{n}}(1 - \bar{r}_{n})}^{+\infty}\omega_{n}^{2}dr \rightarrow +\infty,
\end{equation}
which is a contradiction with (\ref{eqa.15}).

Proof of Theorem \ref{thm2.2} $(iv)$:  We will argue by contradiction, let us assume that there exists $\{u_{n}\}$ such that $u_{n} = u_{\lambda_{n}}$ is the unique solution of (\ref{eq1.1}) in $H_{0,rad}^{1}(D)$ with $\lambda = \lambda_{n} \rightarrow +\infty$ and
\begin{equation} \label{eqa.16}
\liminf_{n \rightarrow +\infty}\int_{D}u_{n}^{2} > 0.
\end{equation}
By $(A_{1})$, $\bar{r}_{n}$ is bounded where $\bar{r}_{n}$ is the unique maximum point of $u_{n}$. When $N = 2$ and $p > 6$, up to a subsequence, (\ref{eqa.12}) shows that for large $n$,
\begin{equation}
\int_{D}u_{n}^{2} \leq 2C\lambda_{n}^{\frac{2}{p-2} - \frac{1}{2}}\bar{r}_{n}^{N-1}\int_{\sqrt{\lambda_{n}}(1 - \bar{r}_{n})}^{+\infty}\omega_{n}^{2}dr \rightarrow 0,
\end{equation}
which is a contradiction with (\ref{eqa.16}).
\qed\vskip 5pt

Our second main result in section reads as follows:

\begin{theorem}
\label{thm2.3} Assume that $(A_{2})$ holds and let $u_{\lambda}$ be the unique solution of (\ref{eq1.1}) in $H_{0,rad}^{1}(D)$ with $\lambda \rightarrow +\infty$. Then $\int_{D}u_{\lambda}^{2} \rightarrow 0$.
\end{theorem}

\begin{proposition}
\label{propB.1} (Pohozaev-type identity) If $u \in H_{0}^{1}(D)$ solves (\ref{eq1.1}), then we have the following integral identity
\begin{eqnarray} \label{eqb.1}
\frac{N-2}{2}\int_{D}|\nabla u|^{2} &-& \frac{N+1}{2}\int_{D}\frac{|\nabla u|^{2}}{|x|} + \int_{D}\frac{|x \cdot \nabla u|^{2}}{|x|^{3}} \nonumber \\
&+& \frac{N\lambda}{2}\int_{D}u^{2}- \frac{(N-1)\lambda}{2}\int_{D}\frac{u^{2}}{|x|} = \frac{N}{p}\int_{D}u^{p} - \frac{N-1}{p}\int_{D}\frac{u^{p}}{|x|}.
\end{eqnarray}
\end{proposition}

\textit{Proof.  } Standard arguments imply that $u \in H_{loc}^{2}(\overline{D})$. It is clear that
$$
\lambda u x\cdot\nabla u = \lambda \mathrm{div}(x\frac{u^{2}}{2}) - \frac{N\lambda}{2}u^{2},
$$
$$
u^{p-1} x\cdot\nabla u = \mathrm{div}(x\frac{u^{p}}{p}) - \frac{N}{p}u^{p},
$$
$$
\Delta u x\cdot\nabla u = \mathrm{div}(\nabla ux\cdot\nabla u - x\frac{|\nabla u|^{2}}{2}) + \frac{N-2}{2}|\nabla u|^{2}.
$$
Integrating by parts $\Delta u - \lambda u + u^{p-1} = 0$ in $T_{1}^{R} = \{x \in \mathbb{R}^{N}: 1 < |x| < R\}$, we obtain
\begin{eqnarray} \label{eqa.17}
&& \frac{N-2}{2}\int_{T_{1}^{R}}|\nabla u|^{2} + \frac{N\lambda}{2}\int_{T_{1}^{R}}u^{2} - \frac{N}{p}\int_{T_{1}^{R}}u^{p} \nonumber \\
&=& \frac{1}{2}\int_{\partial B_{1}}|\nabla u|^{2}d\sigma - \int_{\partial B_{R}}(\frac{1}{R}|\sigma \cdot \nabla u|^{2} - \frac{R}{2}|\nabla u|^{2}  - \frac{R}{2}\lambda u^{2} + \frac{R}{p}u^{p})d\sigma.
\end{eqnarray}
Since $\nabla u \in L^{2}(D)$, $u \in L^{2}(D) \cap L^{p}(D)$,
\begin{equation}
\int_{D}3|\nabla u|^{2}  + \lambda u^{2} + \frac{2}{p}|u|^{p} = \int_{1}^{+\infty}\int_{\partial B_{R}}(3|\nabla u|^{2}  + \lambda u^{2} + \frac{2}{p}|u|^{p})d\sigma dR < +\infty.
\end{equation}
Consequently, there exists a sequence $R_{n} \rightarrow +\infty$ such that
$$
R_{n}\int_{\partial B_{R_{n}}}(3|\nabla u|^{2}  + \lambda u^{2} + \frac{2}{p}|u|^{p})d\sigma \rightarrow 0.
$$
Take $R = R_{n}$ in (\ref{eqa.17}), and let $n \rightarrow +\infty$, we have
\begin{eqnarray} \label{eqa.18}
\frac{N-2}{2}\int_{D}|\nabla u|^{2} + \frac{N\lambda}{2}\int_{D}u^{2} - \frac{N}{p}\int_{D}u^{p} = \frac{1}{2}\int_{\partial B_{1}}|\nabla u|^{2}d\sigma.
\end{eqnarray}

Similarly,
$$
\lambda u \frac{x}{|x|}\cdot\nabla u = \lambda \mathrm{div}(\frac{x}{|x|}\frac{u^{2}}{2}) - \frac{(N-1)\lambda}{2}\frac{u^{2}}{|x|},
$$
$$
u^{p-1} \frac{x}{|x|}\cdot\nabla u = \mathrm{div}(\frac{x}{|x|}\frac{u^{p}}{p}) - \frac{N-1}{p}\frac{u^{p}}{|x|},
$$
$$
\Delta u \frac{x}{|x|}\cdot\nabla u = \mathrm{div}(\nabla u\frac{x}{|x|}\cdot\nabla u - \frac{x}{|x|}\frac{|\nabla u|^{2}}{2}) + \frac{N+1}{2}\frac{|\nabla u|^{2}}{|x|} - \frac{|x \cdot \nabla u|^{2}}{|x|^{3}}.
$$
Integrating by parts $\Delta u - \lambda u + u^{p-1} = 0$ in $T_{1}^{R} = \{x \in \mathbb{R}^{N}: 1 < |x| < R\}$, we obtain
\begin{eqnarray} \label{eqa.19}
&& \frac{N+1}{2}\int_{T_{1}^{R}}\frac{|\nabla u|^{2}}{|x|} - \int_{T_{1}^{R}}\frac{|x \cdot \nabla u|^{2}}{|x|^{3}} + \frac{(N-1)\lambda}{2}\int_{T_{1}^{R}}\frac{u^{2}}{|x|} - \frac{N-1}{p}\int_{T_{1}^{R}}\frac{u^{p}}{|x|} \nonumber \\
&=& \frac{1}{2}\int_{\partial B_{1}}|\nabla u|^{2}d\sigma - \int_{\partial B_{R}}(\frac{1}{R^{2}}|\sigma \cdot \nabla u|^{2} - \frac{1}{2}|\nabla u|^{2}  - \frac{1}{2}\lambda u^{2} + \frac{1}{p}u^{p})d\sigma.
\end{eqnarray}
Take $R = R_{n} \rightarrow +\infty$ in (\ref{eqa.19}) with
$$
\int_{\partial B_{R_{n}}}(3|\nabla u|^{2}  + \lambda u^{2} + \frac{2}{p}|u|^{p})d\sigma \rightarrow 0,
$$
and let $n \rightarrow +\infty$, we have
\begin{eqnarray} \label{eqa.20}
\frac{N+1}{2}\int_{D}\frac{|\nabla u|^{2}}{|x|} &-& \int_{D}\frac{|x \cdot \nabla u|^{2}}{|x|^{3}} \nonumber \\
&+& \frac{(N-1)\lambda}{2}\int_{D}\frac{u^{2}}{|x|} - \frac{N-1}{p}\int_{D}\frac{u^{p}}{|x|} = \frac{1}{2}\int_{\partial B_{1}}|\nabla u|^{2}d\sigma.
\end{eqnarray}
Then (\ref{eqb.1}) follows from (\ref{eqa.18}) and (\ref{eqa.20}). \qed\vskip 5pt

\begin{lemma}
\label{lemB.2} Assume that $(A_{2})$ holds and let $u_{\lambda}$ be the unique solution of (\ref{eq1.1}) in $H_{0,rad}^{1}(D)$ with $\lambda \rightarrow +\infty$. Then for any $\epsilon > 0$, there exists $K > 0$ such that for all $\lambda > K$, we have
\begin{equation} \label{eqb.2}
\int_{D}\frac{|\nabla u_{\lambda}|^{2}}{|x|} \leq \epsilon \int_{D}|\nabla u_{\lambda}|^{2},
\end{equation}
\begin{equation} \label{eqb.3}
\int_{D}\frac{u_{\lambda}^{2}}{|x|} \leq \epsilon \int_{D}u_{\lambda}^{2},
\end{equation}
\begin{equation} \label{eqb.4}
\int_{D}\frac{u_{\lambda}^{p}}{|x|} \leq \epsilon \int_{D}u_{\lambda}^{p}.
\end{equation}
\end{lemma}

\textit{Proof.  } We first prove (\ref{eqb.2}). By $(A_{2})$, $\bar{r}_{\lambda} \rightarrow +\infty$ where $\bar{r}_{\lambda}$ is the unique maximum point of $u_{\lambda}$. We claim that for any $R > 1$, there exists $K > 0$ such that for all $\lambda > K$,
$$
\int_{1 < |x| < R}\frac{|\nabla u_{\lambda}|^{2}}{|x|} \leq \int_{|x| \geq R}\frac{|\nabla u_{\lambda}|^{2}}{|x|}.
$$
Otherwise, there exists some $R > 1$ and $u_{n} = u_{\lambda_{n}}$ with $\lambda_{n} \rightarrow +\infty$ such that
\begin{equation} \label{eqb.6}
\int_{1 < |x| < R}\frac{|\nabla u_{n}|^{2}}{|x|} > \int_{|x| \geq R}\frac{|\nabla u_{n}|^{2}}{|x|}.
\end{equation}
On the other hand,
\begin{eqnarray} \label{eqb.7}
\frac{\int_{1 < |x| < R}\frac{|\nabla u_{n}|^{2}}{|x|}}{\int_{|x| \geq R}\frac{|\nabla u_{n}|^{2}}{|x|}} &=& \frac{\int_{1}^{R}r^{N-2}(u_{n}')^{2}dr}{\int_{R}^{+\infty}r^{N-2}(u_{n}')^{2}dr} \nonumber \\
&\leq& \frac{\int_{1}^{R}(u_{n}')^{2}dr}{\int_{R}^{+\infty}(u_{n}')^{2}dr} \nonumber \\
&=& \frac{\int_{\sqrt{\lambda_{n}}(1 - \bar{r}_{n})}^{\sqrt{\lambda_{n}}(R - \bar{r}_{n})}(\omega_{n}')^{2}dr}{\int_{\sqrt{\lambda_{n}}(R - \bar{r}_{n})}^{+\infty}(\omega_{n}')^{2}dr},
\end{eqnarray}
where $\omega_{n}$ is given by Theorem \ref{thmA.4}. By Theorem \ref{thmA.4}, passing to a subsequence if necessary,
$$
\int_{\sqrt{\lambda_{n}}(R - \bar{r}_{n})}^{+\infty}(\omega_{n}')^{2}dr \rightarrow \int_{-\infty}^{+\infty}(W')^{2}dr > 0
$$
and
$$
\int_{\sqrt{\lambda_{n}}(1 - \bar{r}_{n})}^{\sqrt{\lambda_{n}}(R - \bar{r}_{n})}(\omega_{n}')^{2}dr \rightarrow 0.
$$
Therefore,
$$
\frac{\int_{1 < |x| < R}\frac{|\nabla u_{n}|^{2}}{|x|}}{\int_{|x| \geq R}\frac{|\nabla u_{n}|^{2}}{|x|}} \leq \frac{\int_{\sqrt{\lambda_{n}}(1 - \bar{r}_{n})}^{\sqrt{\lambda_{n}}(R - \bar{r}_{n})}(\omega_{n}')^{2}dr}{\int_{\sqrt{\lambda_{n}}(R - \bar{r}_{n})}^{+\infty}(\omega_{n}')^{2}dr} \rightarrow 0,
$$
which contradicts (\ref{eqb.6}).

For any $\epsilon > 0$, take $R > 2/\epsilon$. There exists $K > 0$ such that for all $\lambda > K$,
\begin{eqnarray} \label{eqb.5}
\int_{D}\frac{|\nabla u_{\lambda}|^{2}}{|x|} &=& \int_{1 < |x| < R}\frac{|\nabla u_{\lambda}|^{2}}{|x|} + \int_{|x| \geq R}\frac{|\nabla u_{\lambda}|^{2}}{|x|} \nonumber \\
&\leq& 2\int_{|x| \geq R}\frac{|\nabla u_{\lambda}|^{2}}{|x|} \nonumber \\
&\leq& \frac{2}{R}\int_{|x| \geq R}|\nabla u_{\lambda}|^{2} \nonumber \\
&\leq& \epsilon\int_{D}|\nabla u_{\lambda}|^{2}.
\end{eqnarray}

Similar arguments yield (\ref{eqb.3}) and (\ref{eqb.4}). \qed\vskip 5pt

\begin{corollary}
\label{corB.3}  Under the hypotheses of Lemma \ref{lemB.2}. For any $\epsilon > 0$, there exists $K > 0$ such that for all $\lambda > K$, we have
\begin{equation} \label{eqb.9}
(N-2 - \epsilon)\int_{D}|\nabla u_{\lambda}|^{2} + (N - \epsilon)\lambda\int_{D}u_{\lambda}^{2} \leq \frac{2N}{p}\int_{D}u_{\lambda}^{p}.
\end{equation}
\end{corollary}

\textit{Proof.  } Corollary \ref{corB.3} can be proved directly by the combination of Proposition \ref{propB.1} and Lemma \ref{lemB.2}.
\qed\vskip 5pt

\textbf{Proof of Theorem \ref{thm2.3}:  } $u_{\lambda}$ solves (\ref{eq1.1}), hence
\begin{equation} \label{eqb.8}
\int_{D}|\nabla u_{\lambda}|^{2} + \lambda\int_{D}u_{\lambda}^{2} = \int_{D}u_{\lambda}^{p}.
\end{equation}
By (\ref{eqb.9}) and (\ref{eqb.8}), we have
$$
2\int_{D}|\nabla u_{\lambda}|^{2} \geq (N-\epsilon-\frac{2N}{p})\int_{D}u_{\lambda}^{p},
$$
$$
2\lambda\int_{D}u_{\lambda}^{2} \leq [\frac{2N}{p}-(N-2-\epsilon)]\int_{D}u_{\lambda}^{p}.
$$
Noting that $u_{\lambda}$ is the ground state of (\ref{eq1.1}), we apply the abstract framework 1 in \cite{Song3} with
$$
\Phi_{\lambda}(u) = \frac{1}{2}\int_{D}(|\nabla u|^{2} + \lambda u^{2}) - \frac{1}{p}\int_{D}u^{p}
$$
and
$$
Q(u) = \frac{1}{2}\int_{D}u^{2}.
$$
Take $\epsilon < \min\{(p-2)N/p, N-(2N+4)/p\}$, then
\begin{eqnarray}
\Phi_{\lambda}(u_{\lambda}) &\geq& \frac{1}{4}(N - \epsilon - \frac{2N + 4}{p})\frac{1}{p}\int_{D}u_{\lambda}^{p} + \frac{\lambda}{2}\int_{D}u_{\lambda}^{2} \nonumber \\
&\geq& (1 + \frac{(N - \epsilon)p - (2N+4)}{2N - (N-2-\epsilon)p})\frac{\lambda}{2}\int_{D}u_{\lambda}^{2},
\end{eqnarray}
i.e. $\Phi_{\lambda}(u_{\lambda}) \geq k\lambda Q(u_{\lambda})$ with $k = 1 + \frac{(N - \epsilon)p - (2N+4)}{2N - (N-2-\epsilon)p} > 1$. By \cite[Theorem 2.2 $(ii)$]{Song3}, $\int_{D}u_{\lambda}^{2} \rightarrow 0$ as $\lambda \rightarrow +\infty$. \qed\vskip 5pt

\begin{remark}
The approach developed in \cite{Song3} can be straightforwardly applied to our case.
\end{remark}

Combining Theorems \ref{thm2.2} and \ref{thm2.3}, we have the following corollary:

\begin{corollary}
\label{corB.4}  Let $u_{\lambda}$ be the unique solution of (\ref{eq1.1}) in $H_{0,rad}^{1}(D)$ with $\lambda \rightarrow +\infty$.

$(i)$ When $N = 2$ and $p < 6$ or $N \geq 3$, $\int_{D}u_{\lambda}^{2} \rightarrow +\infty$.

$(ii)$ When $N = 2$ and $p = 6$, $0 < \liminf_{\lambda \rightarrow +\infty}\int_{D}u_{\lambda}^{2} < +\infty$.

$(iii)$ When $N = 2$ and $p > 6$, $\liminf_{\lambda \rightarrow +\infty}\int_{D}u_{\lambda}^{2} = 0$.
\end{corollary}

\textit{Proof.  } When $N = 2$ and $p < 6$ or $N \geq 3$, Theorem \ref{thm2.2} $(i)$ implies that $\int_{D}u_{\lambda}^{2} \rightarrow +\infty$.

When $N = 2$ and $p = 6$, if $(A_{2})$ holds, Theorem \ref{thm2.2} $(iii)$ implies that $\int_{D}u_{\lambda}^{2} \rightarrow +\infty$. However, by Theorem \ref{thm2.3}, $\int_{D}u_{\lambda}^{2} \rightarrow 0$. Thus we know that $(A_{1})$ holds and then by Theorem \ref{thm2.2} $(ii)$, $0 < \liminf_{\lambda \rightarrow +\infty}\int_{D}u_{\lambda}^{2} < +\infty$.

When $N = 2$ and $p > 6$, if $(A_{1})$ holds, Theorem \ref{thm2.2} $(iv)$ implies that $\liminf_{\lambda \rightarrow +\infty}\int_{D}u_{\lambda}^{2} = 0$. If $(A_{2})$ holds, Theorem \ref{thm2.3} implies that $\int_{D}u_{\lambda}^{2} \rightarrow 0$. Therefore, $\liminf_{\lambda \rightarrow +\infty}\int_{D}u_{\lambda}^{2} = 0$. \qed\vskip 5pt

\section{Asymptotic behaviors when $\lambda \rightarrow 0$}

\begin{lemma}
\label{lemC.1} Let $u_{\lambda}$ be the unique solution of (\ref{eq1.1}) in $H_{0,rad}^{1}(D)$. Then $\lim_{\lambda \rightarrow 0}\|u_{\lambda}\|_{L^{\infty}} \rightarrow 0$.
\end{lemma}

\textit{Proof.  } We can derive the boundedness of $\|u_{\lambda}\|_{L^{\infty}}$ as $\lambda \rightarrow 0$ from similar arguments to Lemma \ref{lemA.2}. Next, we prove that $\lim_{\lambda \rightarrow 0}\|u_{\lambda}\|_{L^{\infty}} \rightarrow 0$. If this is not the case, there exists a sequence $\{u_{n} = u_{\lambda_{n}}\}$ with $\lambda_{n} \rightarrow 0$ such that $\lim_{n \rightarrow +\infty}\|u_{n}\|_{L^{\infty}} \rightarrow M > 0$. Then we consider $\widetilde{u}_{n}(r) = u(r + r_{n})$ where $r_{n} = r_{\lambda_{n}}$.

If $r_{n} \rightarrow +\infty$, we can obtain a solution $\widetilde{u}$ with $\widetilde{u}(0) = M$ satisfying
\begin{equation} \label{eqc.1}
-\widetilde{u}'' = \widetilde{u}^{p-1} \ in \ \mathbb{R}.
\end{equation}
Similar to the proof of Lemma \ref{lemA.3}, we know that (\ref{eqc.1}) admits no solution $\widetilde{u}$ with $\widetilde{u}(0) = M > 0$, which is a contradiction.

If $r_{n} \nrightarrow +\infty$, passing to a subsequence if necessary, we may assume that $r_{n} \rightarrow k \geq 1$. Then we obtain a solution $\widetilde{u}$ with $\widetilde{u}(0) = M$ satisfying
\begin{equation} \label{eqc.2}
-(\widetilde{u}'' + \frac{N-1}{r+k}\widetilde{u}') = \widetilde{u}^{p-1} \ in \ (1-k,+\infty), u(1-k) = 0.
\end{equation}
Let $v(x) = \widetilde{u}(|x|-k)$. Then $v$ satisfies
\begin{equation} \label{eqc.3}
-\Delta v = v^{p-1} \ in \ D, v = 0 \ on \ \partial D.
\end{equation}
By \cite[Proposition 2.1]{GV}, (\ref{eqc.3}) admits no nontrivial positive radial solution, which is a contradiction with $v(x) = M > 0, |x| = k$. \qed\vskip 5pt

\begin{lemma}
\label{lemC.2} Let $u_{\lambda}$ be the unique solution of (\ref{eq1.1}) in $H_{0,rad}^{1}(D)$. Then
\begin{equation}
0 < \liminf_{\lambda \rightarrow 0}\frac{\lambda}{\|u_{\lambda}\|_{L^{\infty}}^{p-2}} \leq \limsup_{\lambda \rightarrow 0}\frac{\lambda}{\|u_{\lambda}\|_{L^{\infty}}^{p-2}} \leq 1.
\end{equation}
\end{lemma}

\textit{Proof.  } Let $\bar{r}_{\lambda}$ be the unique maximum point of $u_{\lambda}$. Then
\begin{equation}
\lambda u_{\lambda}(\bar{r}_{\lambda}) < -u_{\lambda}''(\bar{r}_{\lambda}) + \lambda u_{\lambda}(\bar{r}_{\lambda}) = u_{\lambda}(\bar{r}_{\lambda})^{p-1},
\end{equation}
i.e. $\lambda < u_{\lambda}(\bar{r}_{\lambda})^{p-2}$.

Next we will show that $\liminf_{\lambda \rightarrow 0}\lambda/\|u_{\lambda}\|_{L^{\infty}}^{p-2} > 0$. Suppose on the contrary that there exists $\{\lambda_{n}\}$ with $\lambda_{n} \rightarrow 0$ such that $\lambda_{n}/\|u_{\lambda_{n}}\|_{L^{\infty}}^{p-2} \rightarrow 0$. Similar to the proof of Lemma \ref{lemA.2}, we denote $u_{n} = u_{\lambda_{n}}$ and $\bar{r}_{n} = \bar{r}_{\lambda_{n}}$ and consider
\begin{equation}
v_{n}(r) = \frac{1}{\|u_{n}\|_{L^{\infty}}} u_{n}(\frac{r}{\|u_{n}\|_{L^{\infty}}^{\frac{p-2}{2}}} + \bar{r}_{n}), r \in I_{n} = (\|u_{n}\|_{L^{\infty}}^{\frac{p-2}{2}}(1 - \bar{r}_{n}),+\infty).
\end{equation}

If
$$
\lim_{n \rightarrow +\infty}\|u_{n}\|_{L^{\infty}}^{\frac{p-2}{2}}(1 - \bar{r}_{n}) = -\infty,
$$
we obtain a solution $v$ of
\begin{equation}
-v'' = v^{p-1} \ in \ \mathbb{R}, 0 \leq v \leq 1 \ in \ \mathbb{R},\ v \in C^{2}(\mathbb{R}),\ v(0)=1,
\end{equation}
which is a contradiction with Lemma \ref{lemA.3}.

If
$$
\|u_{n}\|_{L^{\infty}}^{\frac{p-2}{2}}(1 - \bar{r}_{n}) \nrightarrow -\infty,
$$
up to a subsequence, we may assume that
$$
\lim_{n \rightarrow +\infty}\|u_{n}\|_{L^{\infty}}^{\frac{p-2}{2}}(1 - \bar{r}_{n}) = -k \leq 0.
$$
By Lemma \ref{lemC.2},
$$
\lim_{n \rightarrow +\infty}\|u_{n}\|_{L^{\infty}}^{\frac{p-2}{2}}\bar{r}_{n} = k \geq 0.
$$
we obtain a solution $v$ of
\begin{equation}
-(v'' + \frac{N-1}{r + k}v') = v^{p-1} \ in \ (-k,+\infty), v(-k) = 0, 0 \leq v \leq 1, v \in C^{2}([-k,+\infty)),
\end{equation}
with $v(0)=1.$
Noticing that $v(0) = 1$ and $v(-k) = 0$, we know $k > 0$. Let $\widetilde{v}(r) = v(r-k)$. Then $\widetilde{v}$ satisfies
\begin{equation} \label{eqc.10}
-(\widetilde{v}'' + \frac{N-1}{r}\widetilde{v}') = \widetilde{v}^{p-1} \ in \ (0,+\infty), \widetilde{v}(0) = 0.
\end{equation}
If $N \geq 3$, the following Lemma \ref{lemC.7} or \cite[Theorems 5.2, 5.3, 6.2]{Bid} implies that (\ref{eqc.10}) admits no nontrivial positive solution, which is a contradiction with $\widetilde{v}(k) = 1 > 0$. If $N = 2$, Lemma \ref{lemA.1} implies that $\widetilde{v}$ is decreasing in $(k,+\infty)$. However, by \cite[Theorem 7.1]{Bid}, $\widetilde{v}$ oscillates near $+\infty$. We find a contradiction, which completes the proof. \qed\vskip 5pt

\begin{lemma}
\label{lemC.7} If $w \in C^{2}([0,+\infty))$ solves
\begin{equation}
-(w'' + \frac{N-1}{r}w') + \delta w = w^{p-1} \ in \ (0,+\infty), w(0) = 0, w \geq 0,
\end{equation}
with $N > 2$ and $\delta \in \mathbb{R}$. Then $w \equiv 0$.
\end{lemma}

\textit{Proof.  } Suppose on the contrary that $w$ is not always zero. Then $\lim_{r \rightarrow 0}w'(r) = w'(0) > 0$. Since $\delta w(r) - w(r)^{p-1} \rightarrow 0$ as $r \rightarrow 0$, we have $w'' + \frac{N-1}{r}w' \rightarrow 0$. $\forall \epsilon > 0$, $\exists r_{0} = r_{0}(\epsilon) > 0$, $\forall r \in (0,r_{0})$, s.t. $w'' + \frac{N-1}{r}w' \leq \epsilon$. Assume that $w'(r) \geq \alpha > 0$ for all $r \in (0,r_{0})$. We can take $r_{0}$ small enough such that $\alpha > w'(0)/2$ and thus $\epsilon r_{0}/\alpha$ can be chosen small arbitrarily.

Take $\epsilon r_{0}/\alpha < N - 2$ and let $\phi(r) = r^{N - 1 - \epsilon r_{0}/\alpha}w'(r)$. Then for any $r \in (0,r_{0})$,
$$
\phi'(r) = r^{N - 1 - \epsilon r_{0}/\alpha}(w'' + \frac{N-1}{r}w' - \frac{\epsilon r_{0}}{r\alpha}w') \leq 0.
$$
Therefore, $w'(r) \geq Cr^{-(N - 1 - \epsilon r_{0}/\alpha)}$ where $C = r_{0}^{N - 1 - \epsilon r_{0}/\alpha}w'(r_{0}) > 0$. Noticing that
$$
w(r) = \int_{0}^{r}w'(\tau)d\tau \geq C\int_{0}^{r}\tau^{-(N - 1 - \epsilon r_{0}/\alpha)}d\tau, \forall r \in (0,r_{0}),
$$
we derive that $w(r) \rightarrow +\infty$ as $r \rightarrow 0$, contradicting the fact that $w(0) = 0$. \qed\vskip 5pt

\begin{lemma}
\label{lemC.3} Let $u_{\lambda}$ be the unique solution of (\ref{eq1.1}) in $H_{0,rad}^{1}(D)$ and $r_{\lambda}$ be the unique maximum of $u_{\lambda}$. Then $\sqrt{\lambda}r_{\lambda} \rightarrow +\infty$ as $\lambda \rightarrow 0$ if $N \geq 3$ and $r_{\lambda} \rightarrow +\infty$ as $\lambda \rightarrow 0$ if $N = 2$.
\end{lemma}

\textit{Proof.  } First, we assume that $N \geq 3$. By Lemma \ref{lemC.2}, $\sqrt{\lambda} \sim \|u_{\lambda}\|_{L^{\infty}}^{\frac{p-2}{2}}$. Suppose on the contrary that $\|u_{\lambda}\|_{L^{\infty}}^{\frac{p-2}{2}}r_{\lambda} \nrightarrow +\infty$ and there exists a sequence $\{\lambda_{n}\}$ with $\lambda_{n} \rightarrow 0$ such that $\|u_{n}\|_{L^{\infty}}^{\frac{p-2}{2}}r_{n}$ is bounded where $u_{n} = u_{\lambda_{n}}$ and $r_{n} = r_{\lambda_{n}}$. Since $\|u_{n}\|_{L^{\infty}} \rightarrow 0$, we assume that:
$$
\lim_{n \rightarrow +\infty}\|u_{n}\|_{L^{\infty}}^{\frac{p-2}{2}}(\bar{r}_{n}-1) = \lim_{n \rightarrow +\infty}\|u_{n}\|_{L^{\infty}}^{\frac{p-2}{2}}\bar{r}_{n} = k \geq 0.
$$
By Lemma \ref{lemC.2}, passing to a subsequence if necessary, we may assume that $\frac{\lambda}{\|u_{\lambda}\|_{L^{\infty}}^{p-2}} \rightarrow \delta \in (0,1]$. Similar to the proof of Lemma \ref{lemC.2}, we consider
\begin{equation}
v_{n}(r) = \frac{1}{\|u_{n}\|_{L^{\infty}}} u_{n}(\frac{r}{\|u_{n}\|_{L^{\infty}}^{\frac{p-2}{2}}} + \bar{r}_{n}), r \in I_{n} = (\|u_{n}\|_{L^{\infty}}^{\frac{p-2}{2}}(1 - \bar{r}_{n}),+\infty).
\end{equation}
and obtain a solution $v$ with $v(0) = 1$ of
\begin{equation}
-(v'' + \frac{N-1}{r+k}v') + \delta v= v^{p-1} \ in \ (-k,+\infty), v(-k) = 0, 0 \leq v \leq 1, v \in C^{2}([0,+\infty)).
\end{equation}
Noticing that $v(0) = 1$, we have $k > 0$. Let $\widetilde{v}(r) = v(r-k)$. Then $\widetilde{v}$ satisfies
\begin{equation}
-(\widetilde{v}'' + \frac{N-1}{r}\widetilde{v}') + \delta \widetilde{v} = \widetilde{v}^{p-1} \ in \ (0,+\infty), \widetilde{v}(0) = 0,
\end{equation}
contradicting Lemma \ref{lemC.7}.

When $N = 2$, suppose on the contrary that $r_{\lambda} \nrightarrow +\infty$ and there exists a sequence $\{\lambda_{n}\}$ with $\lambda_{n} \rightarrow 0$ such that
$$
\lim_{n \rightarrow +\infty}\|u_{n}\|_{L^{\infty}}^{\frac{p-2}{2}}(\bar{r}_{n}-1) = \lim_{n \rightarrow +\infty}\|u_{n}\|_{L^{\infty}}^{\frac{p-2}{2}}\bar{r}_{n} = 0.
$$
From the above proof, we know that this cannot hold. \qed\vskip 5pt

\begin{theorem}
\label{thmC.4} Let $u_{n}(r) = u_{\lambda_{n}}$ be the unique solution of (\ref{eq1.1}) in $H_{0,rad}^{1}(D)$ with $\lambda = \lambda_{n} \rightarrow 0$ and
\begin{equation}
\omega_{n}(r) = \lambda_{n}^{\frac{1}{2-p}} u_{n}(\frac{r}{\sqrt{\lambda_{n}}} + \bar{r}_{n}), r \in \widetilde{I}_{n} = (\sqrt{\lambda_{n}}(1 - \bar{r}_{n}),+\infty),
\end{equation}
where $\bar{r}_{n} = \bar{r}_{\lambda_{n}}$ is the unique maximum point of $u_{n}$. Then passing to a subsequence if necessary, $\omega_{n}(r) \rightarrow W(r)$ in $H^{1}(\mathbb{R})$ where $W \in H^{1}(\mathbb{R})$ is the unique solution of
\begin{equation}
-W'' + W = W^{p-1} \ in \ \mathbb{R}, W(r) \rightarrow 0 \ as \ |r| \rightarrow +\infty.
\end{equation}
\end{theorem}

\textit{Proof.  } Note that $\sqrt{\lambda_{n}}r_{n} \rightarrow +\infty$ and $\omega_{n}$ satisfies
\begin{equation}
-(\omega_{n}'' + \frac{N-1}{r + \sqrt{\lambda_{n}}\bar{r}_{n}}\omega_{n}') + \omega_{n} = \omega_{n}^{p-1} \ in \ \widetilde{I}_{n}, \omega_{n} = 0 \ on \ \partial \widetilde{I}_{n}.
\end{equation}
Similar to the proof of Theorem \ref{thmA.4}, we can conclude the convergence of $\omega_{n}$. \qed\vskip 5pt

\begin{lemma}
\label{lemC.5} Let $u_{\lambda}$ be the unique solution of (\ref{eq1.1}) in $H_{0,rad}^{1}(D)$ with $\lambda \rightarrow 0$. Then for any $\epsilon > 0$, there exists $\delta > 0$ such that for all $\lambda < \delta$, we have the following identities:
\begin{equation} \label{eqc.4}
\int_{D}\frac{|\nabla u_{\lambda}|^{2}}{|x|} \leq \epsilon \int_{D}|\nabla u_{\lambda}|^{2},
\end{equation}
\begin{equation} \label{eqc.5}
\int_{D}\frac{u_{\lambda}^{2}}{|x|} \leq \epsilon \int_{D}u_{\lambda}^{2},
\end{equation}
\begin{equation} \label{eqc.6}
\int_{D}\frac{u_{\lambda}^{p}}{|x|} \leq \epsilon \int_{D}u_{\lambda}^{p}.
\end{equation}
\end{lemma}

\textit{Proof.  } By Lemma \ref{lemC.3}, $\bar{r}_{\lambda} \rightarrow +\infty$ where $\bar{r}_{\lambda}$ is the unique maximum point of $u_{\lambda}$. Then (\ref{eqc.4}) - (\ref{eqc.6}) can be proved by using similar arguments to Lemma \ref{lemB.2}. \qed\vskip 5pt

\begin{corollary}
\label{corC.6}  Let $u_{\lambda}$ be the unique solution of (\ref{eq1.1}) in $H_{0,rad}^{1}(D)$ with $\lambda \rightarrow 0$. Then for any $\epsilon > 0$, there exists $\delta > 0$ such that for all $\lambda < \delta$, there holds that
\begin{equation} \label{eqc.7}
(N-2 - \epsilon)\int_{D}|\nabla u_{\lambda}|^{2} + (N - \epsilon)\lambda\int_{D}u_{\lambda}^{2} \leq \frac{2N}{p}\int_{D}u_{\lambda}^{p}.
\end{equation}
\begin{equation} \label{eqc.8}
(N-2)\int_{D}|\nabla u_{\lambda}|^{2} + N\lambda\int_{D}u_{\lambda}^{2} \geq \frac{2N}{p}\int_{D}u_{\lambda}^{p}.
\end{equation}
\end{corollary}

\textit{Proof.  } (\ref{eqc.7}) can be proved directly by the combination of Proposition \ref{propB.1} and Lemma \ref{lemC.5}, (\ref{eqc.8}) follows from (\ref{eqa.18}).
\qed\vskip 5pt

\begin{theorem}
\label{thmC.7}  Let $u_{\lambda}$ be the unique solution of (\ref{eq1.1}) in $H_{0,rad}^{1}(D)$ with $\lambda \rightarrow 0$.

$(i)$ If $p < 2 + 4/N$, $\int_{D}u_{\lambda}^{2} \rightarrow 0$.

$(ii)$ If $N \geq 3$, and $p \geq 2 + 4/N$, or $N = 2$, and $p > 4$, $\int_{D}u_{\lambda}^{2} \rightarrow +\infty$.

$(iii)$ If $N = 2$, and $p = 4$,  $\liminf_{\lambda \rightarrow 0}\int_{D}u_{\lambda}^{2} > 0$.
\end{theorem}

\textit{Proof.  } $u_{\lambda}$ solves (\ref{eq1.1}) and thus
\begin{equation} \label{eqc.9}
\int_{D}|\nabla u_{\lambda}|^{2} + \lambda\int_{D}u_{\lambda}^{2} = \int_{D}u_{\lambda}^{p}.
\end{equation}
By (\ref{eqc.7}), (\ref{eqc.8}) and (\ref{eqc.9}), for small $\lambda$, we have
$$
(N-\epsilon-\frac{2N}{p})\int_{D}u_{\lambda}^{p} \leq 2\int_{D}|\nabla u_{\lambda}|^{2} \leq (N-\frac{2N}{p})\int_{D}u_{\lambda}^{p},
$$
$$
[\frac{2N}{p}-(N-2)]\int_{D}u_{\lambda}^{p} \leq 2\lambda\int_{D}u_{\lambda}^{2} \leq [\frac{2N}{p}-(N-2-\epsilon)]\int_{D}u_{\lambda}^{p}.
$$
Noting that $u_{\lambda}$ is the ground state of (\ref{eq1.1}), we apply the abstract framework 1 developed in \cite{Song3} to
$$
\Phi_{\lambda}(u) = \frac{1}{2}\int_{D}(|\nabla u|^{2} + \lambda u^{2}) - \frac{1}{p}\int_{D}u^{p}
$$
and
$$
Q(u) = \frac{1}{2}\int_{D}u^{2}.
$$

Take $\epsilon < \min\{(p-2)N/p, N-(2N+4)/p\}$ when $p > 2 + 4/N$, then
\begin{eqnarray}
\Phi_{\lambda}(u_{\lambda}) &\geq& \frac{1}{4}(N - \epsilon - \frac{2N + 4}{p})\frac{1}{p}\int_{D}u_{\lambda}^{p} + \frac{\lambda}{2}\int_{D}u_{\lambda}^{2} \nonumber \\
&\geq& (1 + \frac{(N - \epsilon)p - (2N+4)}{2N - (N-2-\epsilon)p})\frac{\lambda}{2}\int_{D}u_{\lambda}^{2},
\end{eqnarray}
and
\begin{eqnarray}
\Phi_{\lambda}(u_{\lambda}) &\leq& \frac{1}{4}(N - \frac{2N + 4}{p})\frac{1}{p}\int_{D}u_{\lambda}^{p} + \frac{\lambda}{2}\int_{D}u_{\lambda}^{2} \nonumber \\
&\leq& (1 + \frac{Np - (2N+4)}{2N - (N-2)p})\frac{\lambda}{2}\int_{D}u_{\lambda}^{2},
\end{eqnarray}
i.e. $k\lambda Q(u_{\lambda}) \leq \Phi_{\lambda}(u_{\lambda}) \leq l\lambda Q(u_{\lambda})$ with $k = 1 + \frac{(N - \epsilon)p - (2N+4)}{2N - (N-2-\epsilon)p}$, $l = 1 + \frac{Np - (2N+4)}{2N - (N-2)p}$, $l > k > 1$. By \cite[Theorem 2.2 $(iii)$]{Song3}, $\int_{D}u_{\lambda}^{2} \rightarrow +\infty$ as $\lambda \rightarrow 0$.

Take $\epsilon < 2 - 4/p$ when $p < 2 + 4/N$, then
\begin{eqnarray}
\Phi_{\lambda}(u_{\lambda}) &\geq& \frac{1}{4}(N - \epsilon - \frac{2N + 4}{p})\frac{1}{p}\int_{D}u_{\lambda}^{p} + \frac{\lambda}{2}\int_{D}u_{\lambda}^{2} \nonumber \\
&\geq& (1 + \frac{(N - \epsilon)p - (2N+4)}{2N - (N-2)p})\frac{\lambda}{2}\int_{D}u_{\lambda}^{2},
\end{eqnarray}
and
\begin{eqnarray}
\Phi_{\lambda}(u_{\lambda}) &\leq& \frac{1}{4}(N - \frac{2N + 4}{p})\frac{1}{p}\int_{D}u_{\lambda}^{p} + \frac{\lambda}{2}\int_{D}u_{\lambda}^{2} \nonumber \\
&\leq& (1 + \frac{Np - (2N+4)}{2N - (N-2-\epsilon)p})\frac{\lambda}{2}\int_{D}u_{\lambda}^{2},
\end{eqnarray}
i.e. $k\lambda Q(u_{\lambda}) \leq \Phi_{\lambda}(u_{\lambda}) \leq l\lambda Q(u_{\lambda})$ with $k = 1 + \frac{(N - \epsilon)p - (2N+4)}{2N - (N-2)p}$, $l = 1 + \frac{Np - (2N+4)}{2N - (N-2-\epsilon)p}$, $0< k < l < 1$. By \cite[Theorem 2.2 $(iii)$]{Song3}, $\int_{D}u_{\lambda}^{2} \rightarrow 0$ as $\lambda \rightarrow 0$.

When $p = 2 + 4/N$ and $N \geq 3$, similar to (\ref{eqa.12}),
\begin{eqnarray}
\int_{D}u_{n}^{2} = C\Sigma_{k = 0}^{N-1}C_{N-1}^k(\sqrt{\lambda_{n}}\bar{r}_{n})^{N-1-k}\int_{\sqrt{\lambda_{n}}(1 - \bar{r}_{n})}^{+\infty}\omega_{n}^{2}r^{k}dr,
\end{eqnarray}
for some $C > 0$. By Lemma \ref{lemC.3}, $\int_{D}u_{n}^{2} \rightarrow +\infty$.

When $p = 4$ and $N = 2$, take $\epsilon < 4$. Then
\begin{eqnarray}
\Phi_{\lambda}(u_{\lambda}) &\geq& \frac{- \epsilon}{4p}\int_{D}u_{\lambda}^{p} + \frac{\lambda}{2}\int_{D}u_{\lambda}^{2} \nonumber \\
&\geq& (1 - \frac{\epsilon}{4})\frac{\lambda}{2}\int_{D}u_{\lambda}^{2},
\end{eqnarray}
and
\begin{eqnarray}
\Phi_{\lambda}(u_{\lambda}) \leq \frac{\lambda}{2}\int_{D}u_{\lambda}^{2}
\end{eqnarray}
i.e. $k\lambda Q(u_{\lambda}) \leq \Phi_{\lambda}(u_{\lambda}) \leq \lambda Q(u_{\lambda})$ with $k = 1 - \epsilon/4 \in (0,1)$. Similar to the proof of \cite[Theorem 2.2 $(iii)$]{Song3}, we know that $\liminf_{\lambda \rightarrow 0}\int_{D}u_{\lambda}^{2} > 0$. \qed\vskip 5pt

\textbf{Proof of Theorem \ref{thm1.2}:  } \cite[Theorem 1.1]{FMT} shows that (\ref{eq1.1}) possesses a unique solution $u_{\lambda}$ in $H_{0,rad}^{1}(D)$ for any fixed $\lambda > 0$. Furthermore, $u_{\lambda}$ is the ground state of (\ref{eq1.1}). Then similar to the proof of \cite[Corollary 2.5]{Song}, we know that $\{(\lambda,u_{\lambda}): \lambda > 0\}$ is a continuous curve in $\mathbb{R} \times H_{0,rad}^{1}(D)$. Set $d(\lambda) = \int_{D}u_{\lambda}^{2}$, which is a continuous function of $\lambda > 0$. Note that (\ref{eq1.1}) has a solution with $L^2$ norm $c$ is equvalient to that $d(\lambda) = c$ for some $\lambda > 0$. Also, the number of normalized solutions with mass $c$ is the one of $\{\lambda > 0: d(\lambda) = c\}$.

Proof of $(i)$: If $N \geq 3$, and $2 + 4/N \leq p < 2^{\ast}$  or if $N = 2,$ and $4 < p < 6$, by Corollary \ref{corB.4} $(i)$, $d(\lambda) \rightarrow +\infty$ as $\lambda \rightarrow +\infty$; by Theorem \ref{thmC.7} $(ii)$, $d(\lambda) \rightarrow +\infty$ as $\lambda \rightarrow 0$. Hence, there exists some $\hat{\lambda}$ such that $d(\hat{\lambda}) = \inf_{\lambda > 0}d(\lambda) = \min_{\lambda > 0}d(\lambda) > 0$ since $d(\lambda)$ is continuous. Set
$$
\eta_{1} = d(\hat{\lambda}) = \inf_{\lambda > 0}d(\lambda) > 0.
$$
Then for any $c > \eta_1$, there exist $\lambda > \hat{\lambda}$ and $\widetilde{\lambda} \in (0,\hat{\lambda})$ such that $d(\lambda) = d(\widetilde{\lambda}) = c$,
and for any $c \in (0,\eta_1)$, $d^{-1}(c)$ is empty. This implies the results in Theorem \ref{thm1.2} $(i)$.

Proof of $(ii)$: If $N = 2$, and $p = 6$, by Corollary \ref{corB.4} $(ii)$, $0 < \liminf_{\lambda \rightarrow +\infty}d(\lambda) < +\infty$; by Theorem \ref{thmC.7} $(ii)$, $d(\lambda) \rightarrow +\infty$ as $\lambda \rightarrow 0$. Set
$$
\eta_{2} = \inf_{\lambda > 0}d(\lambda).
$$
If $\inf_{\lambda > 0}d(\lambda) = \liminf_{\lambda \rightarrow +\infty}d(\lambda)$, then $\eta_2 = \liminf_{\lambda \rightarrow +\infty}d(\lambda)$. Otherwise,
$$
\eta_2 = \min_{\lambda > 0}d(\lambda) < \liminf_{\lambda \rightarrow +\infty}d(\lambda).
$$
Hence, for any $c > \eta_2$, $d^{-1}(c)$ is not empty. And for any $c \in (0,\eta_2)$, $d^{-1}(c)$ is empty. Then we can get the results in Theorem \ref{thm1.2} $(ii)$.

Proof of $(iii)$: If $N = 2$, and $p > 6$, by Corollary \ref{corB.4} , $\liminf_{\lambda \rightarrow +\infty}d(\lambda) = 0$; by Theorem \ref{thmC.7} $(ii)$, $d(\lambda) \rightarrow +\infty$ as $\lambda \rightarrow 0$. Hence, for any $c > 0$, $d^{-1}(c)$ is not empty. Then the conclusion of Theorem \ref{thm1.2} $(iii)$ holds true.

Proof of $(iv)$: If $N = 2$, and $p = 4$, by Corollary \ref{corB.4} $(i)$, $d(\lambda) \rightarrow +\infty$ as $\lambda \rightarrow +\infty$; by Theorem \ref{thmC.7} $(iii)$, $\liminf_{\lambda \rightarrow 0}d(\lambda) > 0$. Set
$$
\eta_{3} = \inf_{\lambda > 0}d(\lambda).
$$
If $\inf_{\lambda > 0}d(\lambda) = \liminf_{\lambda \rightarrow 0}d(\lambda)$, then $\eta_3 = \liminf_{\lambda \rightarrow 0}d(\lambda)$. Otherwise,
$$
\eta_3 = \min_{\lambda > 0}d(\lambda) < \liminf_{\lambda \rightarrow 0}d(\lambda).
$$
Hence, for any $c > \eta_3$, $d^{-1}(c)$ is not empty. And for any $c \in (0,\eta_3)$, $d^{-1}(c)$ is empty. Then we can get the results in Theorem \ref{thm1.2} $(iv)$.

Proof of $(v)$: If $2 < p < 2 + 4/N$, by Corollary \ref{corB.4} $(i)$, $d(\lambda) \rightarrow +\infty$ as $\lambda \rightarrow +\infty$; by Theorem \ref{thmC.7} $(i)$, $d(\lambda) \rightarrow 0$ as $\lambda \rightarrow 0$. Hence, for any $c > 0$, $d^{-1}(c)$ is not empty. Then we can get the result in Theorem \ref{thm1.2} $(v)$.
\qed\vskip 5pt

\textbf{Proof of Theorem \ref{thm1.5}:  } Noticing that
$$
\int_{B_{R}^{c}}u_{\lambda,R}^{2} = R^{\frac{4}{2-p}}\int_{B_{R}^{c}}(u_{\lambda R^{2}}(\frac{x}{R}))^{2} = R^{N - \frac{4}{p-2}}\int_{D}(u_{\lambda R^{2}}(x))^{2},
$$
when $2 + 4/N < p < 2^{\ast}$ if $N \geq 3$ and $4 < p < 6$ if $N = 2$, we have
$$
\eta_{1,R} = \inf_{\lambda > 0}\int_{B_{R}^{c}}u_{\lambda,R}^{2} = R^{N - \frac{4}{p-2}}\inf_{\lambda > 0}\int_{D}(u_{\lambda R^{2}}(x))^{2} = R^{N - \frac{4}{p-2}}\eta_{1}.
$$
Since $N - 4/(p-2) > 0$, $\lim_{R \rightarrow 0}\eta_{1,R} = 0$ and $\lim_{R \rightarrow +\infty}\eta_{1,R} = +\infty$.

Similarly, when $N = 2$ and $p = 6$, we have $\lim_{R \rightarrow 0}\eta_{2,R} = 0$ and $\lim_{R \rightarrow +\infty}\eta_{2,R} = +\infty$.

When $p = 2 + 4/N$, $\int_{B_{R}^{c}}u_{\lambda,R}^{2} = \int_{D}(u_{\lambda R^{2}}(x))^{2}$. Therefore, when $N \geq 3$ and $p = 2 + 4/N$, then $\eta_{1,R} \equiv \eta_{1}$; when $N = 2$ and $p = 4$, then $\eta_{3,R} \equiv \eta_{3}$. \qed\vskip 5pt

\section{Orbital stability/instability results}

\begin{lemma}
\label{lemD.1} $\{(\lambda,u_{\lambda}): \lambda > 0\}$ is a $C^{1}$ curve in $\mathbb{R} \times H_{0}^{1}(D)$.
\end{lemma}

\textit{Proof.  } Similar to the proof of \cite[Corollary 2.5]{Song}, we know that $\{(\lambda,u_{\lambda}): \lambda > 0\}$ is a continuous curve in $\mathbb{R} \times H_{0,rad}^{1}(D)$. Furthermore, \cite[Theorem 1.1]{FMT} shows that $u_{\lambda}$ is non-degenerate in $H_{0,rad}^{1}(D)$ for any $\lambda > 0$. Then similar to \cite[Theorem 18]{SS}, we can prove that $\{(\lambda,u_{\lambda}): \lambda > 0\}$ is a $C^{1}$ curve in $\mathbb{R} \times H_{0}^{1}(D)$ by the implicit function Theorem. \qed\vskip 5pt

To study the orbital stability, we lean on the following result, which expresses in our context the abstract theory developed in \cite{GSS}:

\begin{proposition}
\label{propD.2} Assume that $(LWP)$ holds. Then if $\frac{d}{d\lambda}\int_{D}u_{\lambda}^{2} > 0$ (respectively $< 0$), the standing wave $e^{i\lambda t}u_{\lambda}(x)$ is orbitally stable (respectively unstable) in $H_{0,rad}^{1}(D,\mathbb{C})$.
\end{proposition}

\textbf{The proof of Theorem \ref{thm1.3}:  } First, we assume that $2 + 4/N < p < 2^{\ast}$ if $N \geq 3$ and $4 < p < 6$ if $N = 2$ and prove $(i)$. Set $d(\lambda) = \int_{D}u_{\lambda}^{2}$, then $\eta_{1} = \inf_{\lambda > 0}d(\lambda)$. Note that $\lim_{\lambda \rightarrow 0}d(\lambda) = \lim_{\lambda \rightarrow +\infty}d(\lambda) = +\infty$. Take $c \in (\eta_{1},+\infty)$ a regular value of $d$. Assume that $d^{-1}(c) = \{\lambda_{c}^{(1)}, \lambda_{c}^{(2)}, \cdots\}$ with $0< \lambda_{c}^{(1)} < \lambda_{c}^{(2)} < \cdots$, then $d'(\lambda_{c}^{(i)}) \neq 0$, $i = 1, 2, \cdots$. Noticing that $d(\lambda) > d(\lambda_{c}^{(1)})$ in $(0,\lambda_{c}^{(1)})$ and $d(\lambda) < d(\lambda_{c}^{(2)})$ in $(\lambda_{c}^{(1)},\lambda_{c}^{(2)})$, we have $d'(\lambda_{c}^{(1)}) < 0$ and $d'(\lambda_{c}^{(2)}) > 0$. Take $\widetilde{\lambda} = \lambda_{c}^{(1)}$ and $\lambda = \lambda_{c}^{(2)}$. Then according to Proposition \ref{propD.2}, the standing wave $e^{i\lambda t}u_{\lambda}(x)$ is orbitally stable while $e^{i\widetilde{\lambda} t}u_{\widetilde{\lambda}}(x)$ is orbitally unstable in $H_{0,rad}^{1}(D,\mathbb{C})$. By Sard's Theorem, regular values of $d$ are almost every $c \in (\eta_{1},+\infty)$ and thus we complete the proof of $(i)$.

After similar arguments, we can prove $(ii)$ - $(v)$. \qed

\end{document}